\documentclass[10pt,a4paper]{article}
\usepackage[left=3cm,right=3cm]{geometry}
\usepackage[english]{babel}
\usepackage{amssymb}
\usepackage{amsfonts}
\usepackage{amsmath}
\usepackage{color}
\usepackage{amsthm}
\usepackage{graphicx}
\usepackage{upgreek}
\usepackage{multirow}
\usepackage{algorithm}
\usepackage{algorithmic}
\usepackage{tabularx}
\usepackage{stmaryrd} 
\usepackage{tikz}
\usepackage{hyperref}
\usepackage{a4wide}

\newcommand{\be}{\begin{equation}}
\newcommand{\ee}{\end{equation}}
\newcommand{\bi}{\begin{itemize}}
\newcommand{\ei}{\end{itemize}}
\newcommand{\norm}[1]{\left\Vert{#1}\right\Vert} 
\newcommand{\ps}[2]{\left\langle #1,#2\right\rangle}
\newtheorem{proposition}{Proposition}
\newtheorem{definition}{Definition}
\newtheorem{corollaire}{Corollary}
\newtheorem{lemme}{Lemma}

\newtheorem{theoreme}{Theorem}
\newtheorem{remarque}{Remark}

\def\argmin{\textup{argmin}\,}
\newcommand{\prox}{\text{prox}}
\newcommand{\N}{\mathbb{N}}

\newcommand{\R}{\mathbb{R}}

\def\arg{\textup{arg}\,}
\newcolumntype{C}{>{\centering}X}
\newcommand{\off}[1]{}

\begin{document}

\title{
Heavy Ball Momentum for Non-Strongly Convex Optimization}
\author {J.-F. Aujol\footnote{Univ. Bordeaux, Bordeaux INP, CNRS, IMB, UMR 5251, F-33400 Talence, France}
\and 
C. Dossal\footnote{IMT, Univ. Toulouse, INSA Toulouse, Toulouse, France}
\and
H. Labarri\`ere\footnote{MaLGa, DIBRIS, Universit\`a di Genova, Genoa, Italy}
\and 
A. Rondepierre\footnotemark[2] \footnote{LAAS, Univ. Toulouse, CNRS, Toulouse, France}\\
}

\date{\today}

\maketitle



\begin{abstract}

When considering the minimization of a quadratic or strongly convex function, it is well known that first-order methods involving an inertial term weighted by a constant-in-time parameter are particularly efficient (see Polyak \cite{polyak1964some}, Nesterov \cite{nesterov2003introductory}, and references therein). By setting the inertial parameter according to the condition number of the objective function, these methods guarantee a fast exponential decay of the error. We prove that this type of schemes (which are later called Heavy Ball schemes) is relevant in a relaxed setting, i.e. for composite functions satisfying a quadratic growth condition. In particular, we adapt V-FISTA, introduced by Beck in \cite{beck2017first} for strongly convex functions, to this broader class of functions. To the authors' knowledge, the resulting worst-case convergence rates are faster than any other in the literature, including those of FISTA restart schemes. No assumption on the set of minimizers is required and guarantees are also given in the non-optimal case, i.e. when the condition number is not exactly known. This analysis follows the study of the corresponding continuous-time dynamical system (Heavy Ball with friction system), for which new convergence results of the trajectory are shown.


\end{abstract}

\section{Introduction}
In many image processing or statistical problems, the optimization of a convex function $F$ from $\R^N$ to $\R\cup \{+\infty\}$ with a non empty set of minimizers may be needed. In this context, when $N$ is large (i.e. for large scale problems), second order algorithms cannot be used and only gradient or sub-gradient of $F$ can be computed to get a minimizing sequence $(x_n)_{n\in\N}$.

If $F$ is convex, differentiable and has a $L$-Lipschitz gradient, the explicit gradient descent algorithm (GD) with step $s=\frac{1}{L}$ defined by 
\begin{equation*}
x_{n+1}=x_n -s\nabla F(x_n)
\end{equation*}
is a simple first order algorithm that provides a sequence converging to a minimizer $x^*$ of $F$. This method is actually slow on this class of convex functions since its asymptotic convergence rate is
\begin{equation*}
F(x_n)-F(x^*)=\mathcal{O} \left(n^{-1}\right).
\end{equation*} 
This decay rate can be improved when considering $\mu-$strongly convex functions, since the worst-case guarantee is then
\begin{equation}\label{eq:GD_rate_intro}
F(x_n)-F(x^*)=\mathcal{O}\left(e^{-\frac{\mu}{L}n}\right).
\end{equation} 
This asymptotic decay is faster that the one obtained for convex functions but when $\kappa:=\frac{\mu}{L}\ll1$, this decay can still be slow in practice. As $\kappa$ is the inverse of the condition number of $F$, GD is particularly slow for large scale problems.\\
Two remarks can be made about these decays. First, if $F$ is not differentiable but composite, GD can be replaced by the Forward-Backward algorithm and the two decays above are still valid. We provide an exact definition of {\it composite} and {\it Forward-Backward} algorithm in Section 2. Second, the above exponential decay of the error is given under a strong convexity assumption but can be extended under weaker hypotheses such as a quadratic growth condition.

In 1964 Polyak introduces the Heavy Ball (HB) scheme inspired by mechanics, which improves the decay of gradient descent on the class of $C^2$ strongly convex functions by incorporating inertia. This scheme generates a sequence of iterates $(x_n)_{n\in\N}$ ensuring that:
\begin{equation}\label{FGD}
F(x_n)-F(x^*)=\mathcal{O}\left(e^{-4\sqrt{\kappa}n}\right).
\end{equation}
If $\kappa\ll1$, this convergence rate is significantly faster than \eqref{eq:GD_rate_intro} guaranteed by the Forward-Backward algorithm. This theoretical improvement reflects a better performance in practice. At the core of the Polyak's analysis is the fact that in the neighborhood of its unique minimizer, $F$ behaves like a quadratic function. But the $C^2$ assumption is crucial in the Polyak's analysis, and examples of simple $C^1$ strongly convex functions $F$ such that the (HB) provides diverging sequences can be found in \cite{ghadimi2015global}.

In 1983 Nesterov \cite{nesterov1983method} proposes an inertial scheme built to speed up the convergence of GD on the class of convex functions. This acceleration process is at the core of FISTA introduced by Beck and Teboulle \cite{beck2009fast}, which applies to composite functions and provides a sequence such that 
\begin{equation*}
F(x_n)-F(x^*)=\mathcal{O}\left(n^{-2}\right).
\end{equation*}  
The details of this algorithm and convergence rates are given in the Section \ref{sec:def_soa}. The main difference between the Heavy Ball algorithm and the Nesterov scheme is the inertia parameter which is constant over iterations and depends on $\kappa:=\frac{\mu}{L}$ for Heavy Ball while it depends on the iteration number and tends to 1 when $n$ goes to $+\infty$ for FISTA.

Many variations of these schemes have been proposed during the last decade (see Table~\ref{tab:sn_nu} 
for various examples) and the behavior, rates and stability of these various schemes are now well understood. A common approach is to study an associated dynamical system via a Lyapunov analysis before deriving convergence results on the scheme, see e.g. 
\cite{Aujol2023}
\color{black} and the references therein. 

Several Heavy Ball schemes have been proposed to provide fast decays of the type \eqref{FGD} under weaker hypotheses than $C^2$ and strong convexity \cite{nesterov2003introductory,beck2017first,siegel2019accelerated,van2017fastest,taylor2022optimal}. But for all these schemes, a fast geometrical decay such as \eqref{FGD} is achieved only on classes of functions having a unique minimizer: no known inertial scheme achieves such rates on the class of convex functions satisfying a simple quadratic growth condition,
\begin{equation}
\exists \mu>0,~ \forall x\in \R^N,\,\forall x^*\in X^*,\quad F(x)-F(x^*)\geqslant \frac{\mu}{2}d(x,X^*)^2,
\end{equation}   
or equivalently in the convex setting a \L{}ojasiewicz property with parameter $\theta=\frac{1}{2}$, without introducing additional uniqueness hypothesis. In others words, no known inertial scheme provides better asymptotic bounds compared to (GD) within the class of convex functions satisfying a quadratic growth condition. Thus, it remains unclear whether inertia has any real significance for this class of functions, which includes well-known $\ell^1$-regularized functions such as that of \textbf{LASSO problem}. 

The main contribution of this work is to \textbf{provide a Heavy Ball scheme}, similar to Beck's V-FISTA, \textbf{ensuring rates of} $O(e^{-c\sqrt{\kappa}n})$ \textbf{on the class of convex functions satisfying some quadratic growth condition}, where the value of $c$ will be specified later. For a suitable choice of parameters, these worst-case convergence rates are \textbf{faster than that guaranteed by any restarting strategy applied to FISTA}. The method should be parametrized according to the value of $\kappa$, but we prove the exponential decay of the error even if $\kappa$ is overestimated.
    
The paper is organized as follows: in Section \ref{sec:def_soa} we introduce the main geometric assumption made on the function to minimize, namely the quadratic growth condition, and we propose a review of the literature on the convergence rates of inertial algorithms under this condition. Section \ref{sec:contrib} is devoted to our two main theorems proving that Heavy Ball type methods can be properly parameterized to achieve fast exponential decay for this class of function. Section \ref{sec:continu} presents the continuous counterpart of the discrete analysis proposed in Section \ref{sec:contrib} providing a guide to construct the proofs of the theorems presented in Section \ref{sec:contrib}, as well as new results for the convergence rate for the trajectories of the Heavy Ball dynamical system. The proofs have been gathered in Section \ref{sec:proofs}, and the more technical ones are detailed in the Appendix.
 
\section{Geometry of convex functions and inertial algorithms: definitions and state of the art.}\label{sec:def_soa}
 Let us first recall some basic notations and definitions. We assume that $\R^N$ is equipped with the Euclidean scalar product $\langle\cdot,\cdot\rangle$ and the associated norm $\|\cdot\|$. As usual $B(x^*,r)$ denotes the open Euclidean ball with center $x^*\in \R^N$ and radius $r>0$. For any real subset $X\subset \R^N$, the Euclidean distance $d$ is defined as:
$$\forall x\in \R^N,~d(x,X) = \inf_{y\in X} \|x-y\|.$$

 \subsection{Framework and notations}
In this paper we focus on the class of composite functions: $F=f+h$ where $f$ is a convex, differentiable function having a $L$-Lipschitz gradient and $h$ is a proper lower semicontinuous (l.s.c.) convex function whose proximal operator is known. The proximal operator of $h$ is denoted by $\prox_h$ and defined by:
\begin{equation}
    \prox_h(x) = \argmin_{y \in \R^N}{ \left( h(y) + \frac{1}{2} \|y-x \|^2 \right) }.
\end{equation}
For this class of functions a classical minimization algorithm is the Forward-Backward algorithm (FB) whose iterations are described by:
\begin{equation}
x_0\in \R^N,\quad  x_{n+1}=\prox_{sh}(x_n-s \nabla f(x_n)),~s\in \left(0,\frac{2}{L}\right).\label{FB}
\end{equation}
\paragraph{Geometry assumptions.} Without further assumptions on $F$, the convergence decay of the FB algorithm, i.e. the decay of $F(x_n)-F^*$ along the iterates, may be slow. In this paper we are interested in inertial methods, which are among the most effective first order optimization methods, and may ensure a better convergence rates, especially when $F$ is additionally strongly convex:
\begin{definition}[Strong convexity $\mathcal{S}_\mu$] 
Let $F:\R^N\rightarrow \R\cup \{+\infty\}$ be a proper lower semicontinuous convex function. The function $F$ is said $\mu$-strongly convex for some real constant $\mu>0$ if the function $x\mapsto F(x) -\frac{\mu}{2}\|x\|^2$ is convex.
\end{definition}
Weakening this assumption, we consider the class of convex composite functions satisfying some quadratic growth condition:
\begin{definition}[Quadratic growth condition $\mathcal{G}^2_\mu$] 
Let $F:\R^N\rightarrow \R\cup \{+\infty\}$ be a proper lower semicontinuous convex function such that: $X^*=\argmin F\neq \emptyset$ and $F^*=\min F$. The function $F$ satisfies a quadratic growth condition $\mathcal G_\mu^2$ for some real constant $\mu>0$ if:
\begin{equation}
\forall x\in \R^N,~F(x)-F^* \geqslant \frac{\mu}{2}d(x,X^*)^2.\label{Growth:cd}
\end{equation}
\end{definition}
Classically the quadratic growth condition $\mathcal{G}_\mu^2$ can be seen as a relaxation of the strong convexity. Note that satisfying some growth condition does not impose the uniqueness of the minimizer as it does for strong convexity. In the convex setting, the quadratic growth condition $\mathcal{G}_\mu^2$ is equivalent to a global \L ojasiewicz property with an exponent $\frac{1}{2}$ \cite{garrigos2022convergence}. The \L ojasiewicz property \cite{Loja63,Loja93} is a key tool in the mathematical analysis of continuous and discrete dynamical systems, initially introduced to prove the convergence of the trajectories for the gradient flow of analytic functions. An extension to nonsmooth functions has been proposed by Bolte et al. in \cite{Bolte2007loja}:
\begin{definition}
Let $F:\R^N\rightarrow \R\cup\{+\infty\}$ be a lower proper semicontinuous convex function with $X^*=\arg\min~F\neq\emptyset$. The function $F$ has a \L ojasiewicz property of exponent $\theta\in [0,1)$ if for any minimizer $x^*\in X^*$, there exist real constants $c_\ell >0$ and $\varepsilon>0$ such that:
\begin{equation}
    \forall x\in B(x^*,\varepsilon),~\left(F(x)-F(x^*)\right)^\theta \leqslant c_\ell d(0,\partial F(x)).\label{loja}
\end{equation}
The \L ojasiewicz property is said to be global if \eqref{loja} is satisfied for any $x\in \R^N$.
\end{definition}
As shown in \cite[Corollary~9]{bolte2017error}, an example of a non-strongly convex function satisfying the \L{}ojasiewicz property with exponent $\frac{1}{2}$ is the function associated to the LASSO problem
\begin{equation}\label{eq:LASSO}
F(x)=\frac{1}{2}\norm{Ax-y}_2^2+\lambda \norm{x}_1.  
\end{equation}

\paragraph{Inertial methods.} A general inertial optimization method can be described as follows:
\begin{equation}
\forall n \in \N,~\left\{
\begin{gathered}
y_n=x_n+\alpha_n(x_n-x_{n-1}),\\
x_{n+1}=\prox_{sh}\left(y_{n}-s\nabla f(z_{n})\right),
\end{gathered}\right.\label{inertial}
\end{equation}
where $\alpha_n>0$ denotes some friction parameter and $z_n=x_n$ or $y_n$ depending on the considered method. Historically in his seminal work \cite{polyak1964some}, B.T. Polyak proposes a first inertial scheme by choosing a constant friction parameter $\alpha_n=\alpha$ and $z_n=x_n$, for the minimization of $C^2$ strongly convex functions. One of the most popular inertial algorithm is FISTA (for Fast Iterative Shrinkage-Thresholding Algorithm) introduced by Beck and Teboulle in \cite{beck2009fast} to minimize convex composite functions. Inspired by Nesterov's accelerated gradient method proposed in \cite{nesterov1983method}, the friction parameter $\alpha_n$ is defined by:
\begin{equation}
    \alpha_{n} = \frac{t_{n-1}-1}{t_{n}},\qquad z_n = y_n\label{BeckTeboulle}
\end{equation} 
where the sequence $(t_n)_{n\in \N}$ is recursively defined by: $t_0=1$ and $t_{n+1}=\frac{1+\sqrt{1+4t_n^2}}{2}$.
Chambolle and Dossal propose in \cite{chambolle2015convergence} a variant of FISTA defining $\alpha_n=\frac{n-1}{n+\alpha-1}$ for any $n\in\N^*$ where $\alpha\geqslant3$. The original choice of Nesterov can be approximated by setting $\alpha=3$. 

In this paper we consider the family of Heavy Ball algorithms for which the friction parameter $\alpha_n$ is set to a constant $\alpha>0$ and $z_n=y_n$. The term Heavy Ball refers to a family of optimization schemes that can be interpreted as discretizations of the following second-order ordinary differential equation:
\begin{equation}
\ddot x(t) + \alpha \dot x(t) + \nabla F(x(t)) =0,
\end{equation}
which describes the move of a heavy ball in a potential field with a constant friction. The inertia coefficient $\alpha$ has to be parameterized according to the geometry of $F$ to get an optimal convergence rate. For the class of $\mu$-strongly convex functions, Beck in \cite[Chapter~10.7.7]{beck2017first} proposes the following choice (following Nesterov's choice \cite{nesterov2003introductory}):
\begin{equation}
\alpha =\frac{1-\sqrt{\kappa}}{1+\sqrt{\kappa}},\qquad z_n=y_n
\end{equation}
where $\kappa=\frac{\mu}{L}$, leading to the algorithm V-FISTA (seen as a variant of FISTA by the author).


\off{
Beck and Teboulle, based on the ideas of Nesterov's acceleration, propose an accelerated version FISTA (Fast Iterative Shrinkage-Thresholding Algorithm)[13]:
\begin{equation}
 y_n=x_n+ \frac{n}{n+3}(x_n - x_{n-1}),\quad
x_{n+1}=\prox_{sh}(y_n-s \nabla f(y_n)).
\end{equation}
In this paper we consider the variant of FISTA proposed by Chambolle and Dossal in \cite{chambolle2015convergence} and denoted by FISTA:
\begin{equation}
 y_n=x_n+ \frac{n}{n+\alpha}(x_n - x_{n-1}),\quad x_{n+1}=\prox_{sh}(y_n-s \nabla f(y_n))
\end{equation}
that ensures in addition the weak convergence of the iterates (when $\alpha > 3)$. 
}


\subsection{Convergence rate of inertial algorithms under quadratic growth assumptions}
In this section, we give a comprehensive overview of the literature results for first-order methods applied to convex functions satisfying additional growth assumptions. We provide Table \ref{tab:sn_nu} summarizing the theoretical guarantees of Heavy Ball schemes and point out that the uniqueness of the minimizer is required in almost all previous works.

\paragraph{Forward-Backward.} When the function $F$ to minimize satisfies some additional growth assumption $\mathcal{G}^2_\mu$, Garrigos et al. \cite{garrigos2022convergence} prove that the Forward-Backward algorithm \eqref{FB} provides better theoretical guarantees than in the general convex case. More precisely, they show that the function values achieve an exponential decay $F(x_n)-F^*=\mathcal{O}\left(e^{-\frac{\mu}{4L}n}\right)$ along the iterates of Forward-Backward, without any assumption on the set of minimizers $X^*$. 
Observe that this convergence rate depends on the ratio $\kappa=\frac{\mu}{L}$ which represents the inverse of the conditioning of $F$ and can be very small in large-scale optimization. Note also that Necoara et al. proved similar results for the projected gradient algorithm in \cite{necoara2019linear}. 

\paragraph{Nesterov's acceleration.} While Nesterov's acceleration allows for speeding up gradient-based algorithms for the class of convex functions, it is less clear for the class of convex functions satisfying some quadratic growth condition. Considering FISTA, in its historical form by Beck and Teboulle \eqref{BeckTeboulle} or its variant introduced by Chambolle and Dossal in \cite{chambolle2015convergence}, the convergence rate is still polynomial for the class of convex functions satisfying some quadratic growth condition. Note however that considering the variant of FISTA introduced by \cite{chambolle2015convergence}, Aujol et al. prove in \cite{Aujol2023} that the sequence $(x_n)_{n\in\N}$ provided by \eqref{inertial} with $\alpha_n=\frac{n-1}{n+\alpha-1}$ and $\alpha\geqslant3$ satisfies:
\begin{equation}
    F(x_n)-F^*=\mathcal{O}\left(n^{-\frac{2\alpha}{3}}\right).
\end{equation}
which is better than the rate in the convex setting $\mathcal{O}\left(n^{-2}\right)$ from \cite{nesterov1983method,beck2009fast} which is in fact $o\left(n^{-2}\right)$ as proved by Attouch and Peypouquet \cite{attouch2016rate}. Although this decay is not exponential, the authors show that the friction parameter $\alpha$ can be set according to some desired accuracy, and in that case the number of iterations required to achieve this accuracy is comparable to methods ensuring a fast exponential decay of the error, i.e. a exponential decay depending on $\sqrt{\frac{\mu}{L}}$, which can be much faster than an exponential decay rate solely in $\kappa=\frac{\mu}{L}$ as for Forward-Backward. Note that this result originally holds under the assumption that $F$ has a unique minimizer, and was extended without this assumption in \cite{aujol2024strong}.


\paragraph{Restarting strategies.} A way to accelerate the convergence of FISTA for the class of composite convex functions having some quadratic growth property, is to use restart strategies. The idea of this approach is to take benefit of inertia while avoiding oscillations by re-initializing the inertia parameter to zero when some restart condition is verified. Empiric restart rules have been proposed by Giselsson and Boyd \cite{giselsson2014monotonicity} or O'Donoghue and Candès \cite{o2015adaptive}, offering an improved convergence of FISTA in practice but without theoretical guarantees. Elementary computations show that re-initializing the inertia parameter every $\lfloor 2e\sqrt{\frac{L}{\mu}}\rfloor$ iterations allows the resulting sequence to satisfy:
\begin{equation}
    F(x_n)-F^*=\mathcal{O}\left(e^{-\frac{1}{e}\sqrt{\kappa}n}\right).
\end{equation}
This convergence rate is actually the fastest in the literature for restart methods and does not require the uniqueness of the minimizer. But note that it requires knowledge of the value of $\mu$, see e.g. \cite{nesterov2013gradient,o2015adaptive,necoara2019linear}. Recently, adaptive restart schemes have been developed aiming at exploiting the geometry assumption $\mathcal{G}_\mu^2$ to derive improved convergence rates without knowing exactly the growth parameter $\mu$: Fercoq and Qu \cite{fercoq2019adaptive}, Alamo et al. \cite{alamo2022restart}, Aujol et al. \cite{aujol2021restart,aujol2023parameter} introduce restart schemes ensuring a fast exponential decay of the error (i.e. depending on $\sqrt{\kappa})$. The schemes having the best theoretical guarantees in this setting are that proposed by Alamo et al. in \cite{alamo2022restart} ($F(x_n)-F^*=\mathcal{O}\left(e^{-\frac{\ln(15)}{4e}\sqrt{\kappa}n}\right)$) and the method introduced by Renegar and Grimmer in \cite{renegar2022simple} ($F(x_n)-F^*=\mathcal{O}\left(e^{-\frac{1}{2\sqrt{2}}\sqrt{\kappa}n}\right)$). As the optimal periodic restart, no uniqueness assumption is needed on the set of minimizers of $F$ to obtain these guarantees.


\paragraph{Heavy Ball schemes.} In contrast to FISTA and Nesterov's accelerated gradient method, Heavy Ball type schemes are designed for convex functions satisfying additional growth assumptions such as the $\mu$-strong convexity. To this end, these methods require to be calibrated according to the growth parameter $\mu$. In his seminal paper \cite{polyak1964some}, Polyak introduces the first Heavy Ball method for $C^2$ $\mu$-strongly convex functions which guarantees a convergence rate of the error of $\mathcal{O}\left(e^{-4\sqrt{\kappa}n}\right)$. This decay rate is significantly fast but relies strongly on the $C^2$ assumption. Ghadimi et al. in \cite{ghadimi2015global} provide a $C^1$ convex function such that this method does not converge. Nesterov proposes in \cite{nesterov2003introductory} the accelerated gradient method for strongly convex functions which only requires a $C^1$ assumption ensuring that the error decreases as $\mathcal{O}\left(e^{-\sqrt{\kappa}n}\right)$. In this setting, several Heavy Ball schemes have been proposed such as Siegel's Heavy Ball method \cite{siegel2019accelerated} and the geometric descent method \cite{bubeck2015geometric} which have the same theoretical asymptotic guarantees as Nesterov's accelerated gradient method for strongly convex functions, the Heavy Ball method by Aujol et al. \cite{aujol2022convergenceQC} for strongly convex functions (which we will denote ADR-$\mathcal{S}_\mu$ Heavy Ball), the triple momentum method by Van Scoy et al. \cite{van2017fastest} and ITEM by Taylor and Drori \cite{taylor2022optimal}. The latter two schemes are built thanks to the Performance Estimation Problem approach introduced by Drori and Teboulle \cite{Drori2014} and they provide the best bounds on the error for this class of function and first-order methods ($\mathcal{O}\left(e^{-2\sqrt{\kappa}n}\right)$). Some of these schemes can be adapted to composite optimization as detailed in Table \ref{tab:sn_nu}. Note that Beck generalizes Nesterov's accelerated gradient method for strongly convex functions to composite optimization in \cite{beck2017first} with V-FISTA proving the same theoretical convergence rate of the error.

{\renewcommand{\arraystretch}{1.4}
\begin{table}[h]
\footnotesize
\centering
 \begin{tabularx}{\linewidth}{|>{\centering}X|>{\centering}X|>{\centering}X|>{\centering}X|}
  \hline
  \textbf{Algorithm} & \textbf{Reference} & \textbf{Assumption on $F$} & \textbf{Convergence rate of $F(x_n)-F^*$} \tabularnewline
  \hline
  Polyak's Heavy Ball & Polyak \cite{polyak1964some} & $\mathcal{S}_\mu$ and $C^2$ & $\mathcal{O}\left(e^{-4\sqrt{\kappa}n}\right)$\tabularnewline
  \hline
  Nesterov's accelerated gradient method for strongly convex functions & Nesterov \cite{nesterov2003introductory} \\Necoara et al. \cite{necoara2019linear} & $\mathcal{S}_\mu$ and $C^1$\\$\mathcal{Q}_\mu$, $C^1$ and uniqueness of the projection of the iterates onto $X^*$&$\mathcal{O}\left(e^{-\sqrt{\kappa}n}\right)$\tabularnewline 
  \hline
  Geometric descent method & Bubeck et al. \cite{bubeck2015geometric}\\Chen et al. \cite{chen2017geometric} & $\mathcal{S}_\mu$\\Adapted to composite optimization & $\mathcal{O}\left(e^{-\sqrt{\kappa}n}\right)$\tabularnewline
  \hline
  Triple momentum method & Van Scoy et al. \cite{van2017fastest} & $\mathcal{S}_\mu$ and $C^1$ & $\mathcal{O}\left(e^{-2\sqrt{\kappa}n}\right)$ \tabularnewline
  \hline
  ITEM & Taylor, Drori \cite{taylor2022optimal} & $\mathcal{S}_\mu$ and $C^1$ & $\mathcal{O}\left(e^{-2\sqrt{\kappa}n}\right)$ \tabularnewline
  \hline
  Polyak's Heavy Ball with general friction & Ghadimi et al. \cite{ghadimi2015global} & $\mathcal{S}_\mu$ and $C^1$ & $\mathcal{O}\left(e^{-\kappa n}\right)$\tabularnewline
  \hline
  Siegel's Heavy Ball & Siegel \cite{siegel2019accelerated} & $\mathcal{S}_\mu$ and $C^1$\\Adapted to composite optimization& $\mathcal{O}\left(e^{-\sqrt{\kappa}n}\right)$\tabularnewline\hline
  V-FISTA & Beck \cite{beck2017first} & $\mathcal{S}_\mu$\\Adapted to composite optimization & $\mathcal{O}\left(e^{-\sqrt{\kappa}n}\right)$\tabularnewline\hline
  ADR-$\mathcal{S}_\mu$ Heavy Ball & Aujol et al. \cite{aujol2022convergenceQC} & $\mathcal{S}_\mu$\\Adapted to composite optimization & $\mathcal{O}\left(e^{\left(-\sqrt{2\kappa}+\mathcal{O}(\kappa)\right)n}\right)$\tabularnewline\hline
  ADR-$\mathcal{G}^2_\mu$ Heavy Ball & Aujol et al. \cite{aujol2022convergenceL} & $\mathcal{G}^2_\mu$ and uniqueness of the minimizer\\Adapted to composite optimization & \tiny{$\mathcal{O}\left(e^{(-(2-\sqrt{2})\sqrt{\kappa}+\mathcal{O}(\kappa))n}\right)$}\tabularnewline\hline
  ADR-$\mathcal{G}^2_\mu$ Heavy Ball & Aujol et al. \cite{aujol2022convergenceL} & $\mathcal{G}^2_\mu$\\Adapted to composite optimization & $\mathcal{O}\left(e^{(-\kappa+\varepsilon+\mathcal{O}(\kappa))n}\right)$\tabularnewline\hline
\end{tabularx}\normalsize
\caption{Convergence rate of $F(x_n)-F^*$ for Heavy Ball type schemes with various geometry assumptions on $F$.}\label{tab:sn_nu}
\end{table}}

Recently, Heavy Ball type schemes have been studied under weaker geometry assumptions than strong convexity. Necoara et al. prove in \cite{necoara2019linear} that the convergence rate of Nesterov's accelerated gradient method for strongly convex method is actually valid for $C^1$ $\mu$-quasi-strongly convex functions i.e. for functions satisfying:
\begin{equation}
\forall x\in \R^n,  \langle\nabla F(x),x-x^*\rangle \geqslant F(x)-F(x^*) + \frac{\mu}{2}\|x-x^*\|^2,
\color{black}\end{equation}
where $x^*$ denotes the projection onto $X^*$, provided that the iterates share the same projection onto the set of minimizers. In \cite{aujol2022convergenceL}, Aujol et al. build a Heavy Ball type scheme (ADR-$\mathcal{G}^2_\mu$ Heavy Ball) for functions satisfying the quadratic growth assumption $\mathcal{G}^2_\mu$ guaranteeing that
\begin{equation}
    F(x_n)-F^*=\mathcal{O}\left(e^{(-(2-\sqrt{2})\sqrt{\kappa}+\mathcal{O}(\kappa))n}\right),
\end{equation}
as long as $F$ has a unique minimizer. 

Thus, the theoretical guarantees of Heavy Ball type schemes are the best in the literature among first-order methods for functions satisfying growth conditions but they do not hold without assuming the uniqueness of the minimizer. If this hypothesis is not verified, the theoretical convergence rates are similar to those of Forward-Backward, and the relevance of applying such algorithms in this context can therefore be questioned.

\section{Fast exponential decay for Heavy Ball type algorithms}\label{sec:contrib}
We now consider Heavy Ball type methods that can be generically described as variants of the V-FISTA algorithm proposed by Beck in \cite{beck2017first}:
\begin{equation}
\forall n \in \N,~\left\{
\begin{gathered}
x_{n+1}=\prox_{sh}\left(y_{n}-s\nabla f(y_{n})\right),\\
y_{n+1}=x_{n+1}+\alpha(x_{n+1}-x_{n}),
\end{gathered}\right.\label{eq:V-FISTA} \tag{V-FISTA}
\end{equation}
with $x_0\in\R^N$, $s= \frac{1}{L}$, $y_0=x_0$ and any $\alpha>0$. Recall that in the original definition of V-FISTA \cite{beck2017first}, the damping parameter $\alpha$ is set to: $\frac{1-\sqrt{\kappa}}{1+\sqrt{\kappa}}$ where $\kappa=\frac{\mu}{L}$ denotes the inverse of the conditioning of the function $F$ to minimize.

The main contribution in this section is to prove that Heavy Ball methods like \eqref{eq:V-FISTA} can be properly parameterized to achieve fast exponential decay rates (i.e. depending on $\sqrt{\kappa}$) for the class of convex composite functions satisfying some quadratic growth property $\mathcal G^2_\mu$ and without assuming the uniqueness of the minimizer. To the best of our knowledge, this is the first result proving that an inertial method can actually improve the convergence rate of the Forward-Backward algorithm (which is in $O(e^{-\frac{\kappa}{4} n})$) in this setting. In large scale dimension, the inverse $\kappa$ of the conditioning of the function to minimize could be very small, so that decaying in $\kappa$ could be much slower than in $\sqrt{\kappa}$.

We provide two theorems stating worst-case convergence bounds. In Theorem \ref{thm:V-FISTA}, we give guarantees for a well-chosen value of $\alpha$ depending on $\kappa$. More precisely, it is shown that for such a parametrization the error decreases with a fast exponential decay, ensuring the best rate in the literature for this assumption. Theorem \ref{thm:vfista_v2} states more complete results as it provides convergence bounds parametrized by the choice $\alpha$. In particular, this theorem hides improved rates as highlighted in Corollary \ref{cor:tab0}, and gives guarantees when $\alpha$ is poorly chosen (see Corollary \ref{cor:vfista_slow}). The proof of this second theorem is inspired by the proof of Theorem \ref{thm:continu} which provides results on the solution of the Heavy Ball dynamical system. 

\paragraph{Fast exponential decay of the error.}
\begin{theoreme}
\label{thm:V-FISTA}
Let $F=f+h$ where $f$ is a convex differentiable function having a $L$-Lipschitz gradient for some $L>0$, and $h$ a proper convex l.s.c. function. Assume that $F$ satisfies a quadratic growth condition $\mathcal G_\mu^2$ for some real parameter $\mu>0$. Let $(x_n)_{n\in\N}$ be the sequence provided by \eqref{eq:V-FISTA} with $\alpha=1-\frac{5}{3\sqrt{3}}\sqrt{\kappa}$ and $s=\frac{1}{L}$. If $\kappa\leqslant\frac{1}{3}$, then:
\begin{equation}\label{eq:thm1_error}
\forall n\in \N,~F(x_n)-F^*\leqslant \frac{4}{3}\left(1-\frac{2}{3\sqrt{3}}\sqrt{\kappa}\right)^n(F(x_0)-F^*),
\end{equation}
and
\begin{equation}
\|x_n-x_{n-1}\|=\mathcal{O}\left(e^{-\frac{\sqrt{\kappa}}{3\sqrt{3}}n}\right).
\end{equation}
\end{theoreme}

Theorem \ref{thm:V-FISTA}, whose proof can be found in Section~\ref{sec:proof_V-FISTA}, ensures that for a well-chosen parameter $\alpha$ which depends on $\kappa$, the decay of the error along the iterates of \eqref{eq:V-FISTA} is at worst of order $\mathcal{O}\left(e^{-\frac{2}{3\sqrt{3}}\sqrt{\kappa}n}\right)$. Looking back at the results in the literature, this convergence rate is slower than those of most other Heavy Ball schemes. However, we recall that the required assumptions on $F$ in these works (summarized in Table \ref{tab:sn_nu}) are stronger than those needed in Theorem \ref{thm:V-FISTA}. The only method proposed for functions satisfying $\mathcal{G}^2_\mu$, i.e. ADR-$\mathcal{G}^2_\mu$ Heavy Ball \cite{aujol2022convergenceL}, ensures a fast exponential decay of the error: $$F(x_n)-F^*=\mathcal{O}\left(e^{(-(2-\sqrt{2})\sqrt{\kappa}+\mathcal{O}(\kappa))n}\right),$$
if the function $F$ has a unique minimizer. This theoretical decay is faster than \eqref{eq:thm1_error}, but it does not hold if $F$ has multiple minimizers. To the authors' knowledge, the fast exponential decay of the error given by Theorem \ref{thm:V-FISTA} is the first in the literature for Heavy Ball methods in this setting and without any uniqueness assumption on the set of minimizers.

In addition, the guarantee on the decay of the error given in \eqref{eq:thm1_error} is faster than that given by FISTA restarted periodically and optimally as it only ensures (even with some oracle \cite{necoara2019linear})
$$F(x_n)-F^*=\mathcal{O}\left(e^{-\frac{1}{e}\sqrt{\kappa}n}\right).$$
This means that in this setting, \eqref{eq:V-FISTA} is theoretically more relevant than a periodic restart of FISTA when the growth parameter $\mu$ is known. In other words, one should define a constant inertia parameter depending on $\mu$ and $L$ instead of setting an increasing inertia parameter and re-initializing it optimally. 

The second claim of Theorem \ref{thm:V-FISTA} gives an asymptotic control on $\|x_n-x_{n-1}\|$ which ensures that $\sum_{n\in\N^*}\|x_n-x_{n-1}\|<+\infty$. As a consequence, the length of the trajectory of the sequence $\left(x_n\right)_{n\in\N}$ is finite and it converges strongly to a minimizer of the function $F$.

\paragraph{Parametric convergence rates.}
Below is a second theorem about the iterates of \eqref{eq:V-FISTA} which gives stronger and more general results than Theorem \ref{thm:V-FISTA}. The proof is built using the parallel with dynamical systems (see Section \ref{sec:continu}) and is located in Section \ref{sec:proof_v-fista2}.

\begin{theoreme}
\label{thm:vfista_v2}
Let $F=f+h$ where $f$ is a convex differentiable function having a $L$-Lipschitz gradient for some $L>0$, and $h$ a proper convex l.s.c. function. Assume that $F$ satisfies a quadratic growth condition $\mathcal G_\mu^2$ for some real parameter $\mu>0$. Let $(x_n)_{n\in\N}$ be the sequence provided by \eqref{eq:V-FISTA} with $\alpha=1-\omega\sqrt{\kappa}$ where $\omega\in\left(0,\frac{1}{\sqrt{\kappa}}\right)$. Then, for any $n\in\N$:
\begin{equation}
F(x_n)-F^*\leqslant \left(1+(\omega-\tau)^2+(\omega-\tau)\omega\tau\sqrt{\kappa}\right)\left(1-\tau\sqrt{\kappa}+\tau^2\kappa\right)^n(F(x_0)-F^*),\label{eq:thm_v2_1}
\end{equation}
and 
\begin{equation}\label{eq:thm_v2_2}
    \|x_n-x_{n-1}\|=\mathcal{O}\left(e^{-\frac{1}{2}\tau\sqrt{\kappa}\left(1-\tau\sqrt{\kappa}\right)n}\right),
\end{equation}
where $\tau>0$ satisfies the following inequality:
\begin{equation}
\left(1-\omega\sqrt{\kappa}\right)\tau^3-\omega\left(2-\omega\sqrt{\kappa}\right)\tau^2+(\omega^2+2)\tau-\omega\leqslant0.\label{eq:cond_vfista}
\end{equation}
\end{theoreme}

The statements of Theorem \ref{thm:vfista_v2} are less readable than those of Theorem \ref{thm:V-FISTA} but they are actually stronger. The inequality \eqref{eq:cond_vfista} hides the convergence rates which can be obtained for a given $\omega\in\left(0,\frac{1}{\sqrt{\kappa}}\right)$. Observe that the larger $\tau$, the better the convergence rate. The best rates are obtained when:
\begin{equation}\label{eq:cond_vfista_eg}
\left(1-\omega\sqrt{\kappa}\right)\tau^3-\omega\left(2-\omega\sqrt{\kappa}\right)\tau^2+(\omega^2+2)\tau-\omega=0.
\end{equation}
The admissible maximum value of $\tau$ can be thus obtained studying the limit case when $\kappa=0$:
\begin{equation}
\tau^3 - 2 \omega \tau^2 +(\omega^2+2) \tau -\omega =0.
\end{equation}
i.e. (since $\tau >0)$:
\begin{equation}
\omega^2 - \frac{1+2 \tau^2}{\tau} \omega +2 + \tau^2=0.
\end{equation}
This can be seen as a quadratic polynomial in $\omega$, whose discriminant is:
\begin{eqnarray*}
    \Delta &=& \frac{1}{\tau^2} \left(
1 + 4 \tau^4 +4 \tau^2 - 8 \tau^2 - 4 \tau^4
    \right)\\
    &=& \frac{1}{\tau^2} \left(
1 -4 \tau^2
    \right)
    = \frac{1}{\tau^2} \left(
1 -2\tau
    \right)
    (1+2 \tau)
\end{eqnarray*}
Hence the largest value of $\tau$ for which the discriminant satisfies $\Delta\geqslant 0$ is $\tau= \frac{1}{2}$, which corresponds to a limit maximum value of $\omega$ equal to $\omega=\frac{3}{2}$. These two observations highlight that the convergence rates given by Theorem~\ref{thm:vfista_v2}, are faster than that given in Theorem~\ref{thm:V-FISTA} for suitable choices of $\alpha$, as expressed in the following corollary. Note that the convergence guarantees and best choices of $\alpha$ depend on the value of the conditioning number since $\kappa$ appears in Equation \eqref{eq:cond_vfista_eg}.


\begin{corollaire}\label{cor:tab0}
Let $F=f+h$ where $f$ is a convex differentiable function having a $L$-Lipschitz gradient for some $L>0$, and $h$ a proper convex l.s.c. function. Assume that $F$ satisfies a quadratic growth condition $\mathcal G_\mu^2$ for some real parameter $\mu>0$. Let $\kappa=\frac{\mu}{L}$.

Let $(\omega,\tau)\in (\R_+)^2$ be two real parameters chosen such that:
\begin{equation}
\left(1-\omega\sqrt{\kappa}\right)\tau^3-\omega\left(2-\omega\sqrt{\kappa}\right)\tau^2+(\omega^2+2)\tau-\omega=0
\end{equation}
and the value of $\tau$ is maximum. Let $(x_n)_{n\in\N}$ be the sequence provided by \eqref{eq:V-FISTA} with $\alpha=1-\omega\sqrt{\kappa}$. Then:
\begin{equation}
\begin{aligned}
F(x_n)-F^*&\leqslant C\left(1-\sigma\sqrt{\kappa}\right)^n(F(x_0)-F^*)
\end{aligned}
\end{equation}
and
\begin{equation}
    \|x_n-x_{n-1}\|=\mathcal{O}\left(e^{-\frac{\sigma}{2}\sqrt{\kappa}n}\right),
\end{equation}
where: 
\begin{equation*}
    C = 1+(\omega-\tau)^2+(\omega-\tau)\omega\tau\sqrt{\kappa}, \quad     \sigma = \tau -\tau^2\sqrt\kappa,
\end{equation*} 
and there exist three real-valued functions $\varepsilon_i$, $i=1,2,3$  such that: $\lim_{t\to 0}\varepsilon_i(t) =0$ and:
\begin{equation*}
\omega =  \frac{3}{2}-\varepsilon_1(\kappa),\qquad
\tau = \frac{1}{2} -\varepsilon_2(\kappa),\qquad \sigma = \frac{1}{2} -\varepsilon_2(\kappa).
\end{equation*}
\end{corollaire}

Table \ref{tab:rates} provides admissible sets of parameters $(\omega,\tau)$ for Corollary~\ref{cor:tab0}.
\begin{table}[H]
\label{tab:rates}
\centering
\begin{tabular}{||c c c c c||} 
 \hline
 $\kappa$ & $\omega$ & $\tau$ & $\sigma$ & $C$  \\ [0.5ex] 
 \hline\hline
 $1$ & $1.2$ & $\mathbf{0.39}$ & $0.23$ & $2.04$ \\ 
 $\frac{1}{3}$ & $1.32$ & $\mathbf{0.42}$ & $0.31$ & $2.1$ \\
 $10^{-1}$ & $1.39$ & $\mathbf{0.45}$ & $0.38$ & $2.07$ \\
 $10^{-2}$ & $1.46$ & $\mathbf{0.48}$ & $0.45$ & $2.03$ \\
 $10^{-3}$ & $1.49$ & $\mathbf{0.494}$ & $0.486$ & $2.02$ \\
 $10^{-4}$ & $1.495$ & $\mathbf{0.498}$ & $0.495$ & $2.002$ \\[1ex]
 \hline
\end{tabular}
\caption{Admissible sets of parameters for Corollary~\ref{cor:tab0}.}
\end{table}
Thus, Corollary \ref{cor:tab0} provides better convergence rates than Theorem \ref{thm:V-FISTA} (since $\frac{2}{3\sqrt{3}}\approx0.38$). We can remark that the guarantees given by Aujol et al. in \cite{aujol2022convergenceL} for ADR-$\mathcal{G}^2_\mu$ are still better with the additional  assumption that $F$ has a unique minimizer.

\off{One gets Corollary~\ref{cor:tab0} from Corrollary~\ref{cor:tab1} by using equation \eqref{eq:cond_vfista} with $\kappa=0$. Indeed the limit case is when \eqref{eq:cond_vfista} is an equality, and with 
$\kappa=0$ one gets:
\begin{equation}
\tau^3 - 2 \omega \tau^2 +(\omega^2+2) \tau -\omega =0.
\end{equation}
i.e. (since $\tau >0)$:
\begin{equation}
\omega^2 - \frac{1+2 \tau^2}{\tau} \omega +2 + \tau^2=0.
\end{equation}
This can be seen as a quadratic polynomial in $\omega$, whose discriminant is:
\begin{equation}
    \Delta = \frac{1}{\tau^2} \left(
1 + 4 \tau^4 +4 \tau^2 - 8 \tau^2 - 4 \tau^4
    \right)
\end{equation}
i.e.:
\begin{equation}
    \Delta = \frac{1}{\tau^2} \left(
1 -4 \tau^2
    \right)
    = \frac{1}{\tau^2} \left(
1 -2\tau
    \right)
    (1+2 \tau)
\end{equation}
Hence the largest value of $\tau$ for which the discriminant satisfies $\Delta\leq 0$ is $\tau= \frac{1}{2}$, and thus $\omega=\frac{3}{2}$, as given in Corollary~\ref{cor:tab0}.}

\paragraph{Guarantees in suboptimal cases.}

\begin{figure}[ht]
\centering
\includegraphics[width=\textwidth]{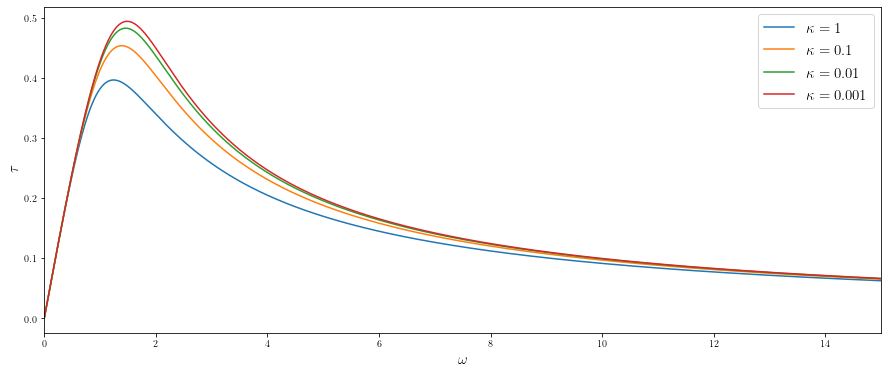}
\caption{Evolution of $\tau$ with respect to $\omega$ for several values of $\kappa$ such that $(\tau,\omega,\kappa)$ satisfy \eqref{eq:cond_vfista_eg}.}
\label{fig:VFISTA_rates}
\end{figure}

In fact, Theorem \ref{thm:vfista_v2} and inequality \eqref{eq:cond_vfista} hide more than improved convergence rates. To illustrate this, we provide a graph displaying the evolution of $\tau$ with respect to $\omega$ and $\kappa$ such that $(\tau,\omega,\kappa)$ satisfy \eqref{eq:cond_vfista_eg} in Figure \ref{fig:VFISTA_rates}. An interactive graph can be found on the link \url{https://www.desmos.com/calculator/syrtiatos6}. We can see on this graph that inequality \eqref{eq:cond_vfista} allows to obtain convergence guarantees even for non-optimal choices of $\alpha$, i.e. large values of $\omega$.

By exploiting this observation, the following corollary provides convergence rates for \eqref{eq:V-FISTA} if $\alpha$ is too small which can be the case if $\mu$ is overestimated. A brief proof is given in Section \ref{sec:proof_cor2}.

\begin{corollaire}\label{cor:vfista_slow}
Let $F=f+h$ where $f$ is a convex differentiable function having a $L$-Lipschitz gradient for some $L>0$, and $h$ a proper convex l.s.c. function. Assume that $F$ satisfies a quadratic growth condition $\mathcal G_\mu^2$ for some real parameter $\mu>0$. Let $(x_n)_{n\in\N}$ be the sequence provided by \eqref{eq:V-FISTA} with $s=\frac{1}{L}$ and $\alpha=1-\theta$ for some $\theta\in\left[\frac{3}{2}\sqrt{\kappa},1\right)$. Then, if $\kappa\leqslant \frac{1}{10}$, 
\begin{equation}
F(x_n)-F^*=\mathcal{O}\left(e^{-\tau\kappa n}\right),
\end{equation}
and 
\begin{equation}
    \|x_n-x_{n-1}\|=\mathcal{O}\left(e^{-\frac{\tau}{2}\kappa n}\right)
\end{equation}
where $\tau=\frac{2}{3\theta}\left(1-\frac{2}{3\theta}\sqrt{\kappa}\right)$.
\end{corollaire}

Corollary \ref{cor:vfista_slow} allows us to derive convergence rates for \eqref{eq:V-FISTA} even if $\alpha$ is far from its optimal value. 
Let us describe two examples:
\begin{itemize}
    \item Suppose that $\alpha=1-C\sqrt{\kappa}$ where $C\geqslant \frac{3}{2}$. Then, applying Corollary \ref{cor:vfista_slow} with $\theta=C\sqrt{\kappa}$, we get that the iterates of \eqref{eq:V-FISTA} for this inertia parameter ensure a decrease of the error in $\mathcal{O}\left(e^{-\tau\sqrt{\kappa} n}\right)$ where $\tau=\frac{2(3C-2)}{9C^2}$. Thus, if we choose $\alpha=1-\frac{3}{2}\sqrt{\tilde\kappa}$ where $\tilde\kappa$ is an upper estimation of $\kappa$, then we get a theoretical guarantee on the error with $C=\frac{3}{2}\sqrt{\frac{\tilde\kappa}{\kappa}}$. In this way, we obtain that if $\tilde\kappa=10\kappa$, then the error decreases as $\mathcal{O}\left(e^{-\tau\sqrt{\kappa}n}\right)$ where $\tau\approx0.12$.
    \item Assume now that $\alpha$ is arbitrarily set to $\alpha=\frac{9}{10}$ without knowing the actual value of $\kappa$. Consequently, we have that $\alpha=1-\theta$ where $\theta=\frac{1}{10}$. If $\theta\geqslant \frac{3}{2}\sqrt{\kappa}$ (i.e. if $\kappa<\frac{1}{225}$), then Corollary \ref{cor:vfista_slow} states that the error along the iterates of \eqref{eq:V-FISTA} for this choice of $\alpha$ decreases in the worst case as $e^{-\tau\kappa n}$ where $\tau=\frac{20}{3}\left(1-\frac{20}{3}\sqrt{\kappa}\right)$ which can be upper bounded by $\frac{100}{27}\approx3.7$. As a consequence, if $F$ is sufficiently ill-conditioned, \eqref{eq:V-FISTA} has better theoretical guarantees for this choice of $\alpha$ than Forward-Backward which ensures that $F(x_n)-F^*=\mathcal{O}\left(e^{-\frac{\kappa}{4} n}\right)$ \cite{garrigos2022convergence}.
\end{itemize}

The convergence guarantees obtained in non-optimal cases and the robustness to an overestimation of the growth parameter are the main contributions of this work. To the authors' knowledge, there are no such results in the literature. This provides a better understanding of the behavior of the iterates of \eqref{eq:V-FISTA} for a wide range of values of $\alpha$.

%
%


\section{Convergence rates for the trajectories of the Heavy Ball dynamical system}\label{sec:continu}

The so called Heavy Ball equation 
\begin{equation}\label{HB_ODE}
\ddot{x}(t)+\alpha \dot{x}(t)+\nabla F(x(t))=0
\end{equation}
where $\alpha>0$ denotes the friction parameter, has been studied from decades now. See for instance Attouch et al \cite{attouch2000heavy} for a general study of the dynamical system : existence of solutions, link with the mechanical system and convergence of the trajectories to critical points if no strong assumptions are made on $F$. The main result of this section is to prove that if $F$ is convex, differentiable and satisfies a quadratic growth condition, the solution of \eqref{HB_ODE} ensures a fast exponential decay. A crucial point here, is that the uniqueness of the minimizer of $F$ is not needed to get such fast rates. 
Before stating Theorem \ref{thm:continu}, we give an overview of the literature. We highlight that slower exponential decays are already known in this setting, that a fast decay is known if the minimizer of $F$ is supposed to be unique and that various other decays can be achieved in non convex settings. Note that the proof of Theorem \ref{thm:continu} was the guideline for the analysis that led to Theorem \ref{thm:vfista_v2}. 

\subsection{State of the art}

The first study of the convergence rate of $F(x(t))-F^*$, where $x$ is a solution of \eqref{HB_ODE}, under strong convexity analysis is due to Polyak \cite{polyak1964some} . In his seminal work he observes that if the function $F$ is a quadratic function, $F(x)=\norm{Ax}^2$, the solution of \eqref{HB_ODE} ensures the fastest decay of $F(x(t))-F^*$ when $\alpha=\sqrt{\mu}$ where $\mu$ is the smallest non negative eigenvalue of $A^\top A$. He deduced that if $F$ is $C^2$ and $\mu$-strongly convex $(\mathcal{S}_{\mu})$, the solution of \eqref{eq:HB_ODE} satisfies 
\begin{equation}
    F(x(t))-F(x^*)=\mathcal{O}(e^{-2\sqrt{\mu}t}).
\end{equation}
Both $C^2$ and strong convexity are necessary to achieve such a decay. During the last decades several works provide various bounds depending on geometrical assumptions made on $F$. A summary is given in Table \ref{tab:hb_nu}.

If $F$ is $\mu$-strongly convex and $C^1$ it was first proved that the decay was $\mathcal{O}(e^{-\sqrt{\mu}t})$, see for example Siegel \cite{siegel2019accelerated} for a simple proof. Aujol et al. \cite{aujol2022convergenceQC} extend this former result giving a better rate for functions that are quasi-strongly convex and have a unique minimizer, and a weaker convergence rate if $F$ satisfies only a quadratic growth condition and has a unique minimizer, see Table \ref{tab:hb_nu} for more details. 
All these results ensure fast exponential decays of $F(x(t))-F^*$ and assume the convexity of $F$, a quadratic growth condition and a uniqueness of the minimizer of $F$.

In \cite{begout2015damped}, Bégout et al. provide several results on the trajectory  solution of \eqref{eq:HB_ODE} if $F$ is a $C^2$ function satisfying some \L{}ojasiewicz property, but is not necessarily convex. The authors prove that the trajectory strongly converges to a minimizer of $F$, and provide several decay rates. In particular, if $F$ satisfies a \L{}ojasiewicz property with an exponent $\theta=\frac{1}{2}$, the error along the trajectory decreases exponentially. More recent works by Polyak and Shcherbakov \cite{polyakshch}, Apidopoulos et al \cite{apidopoulos2022convergence} and Kassing and Weissmann \cite{kassing2024polyak} also prove the exponential decay of the error under the \L{}ojasiewicz assumption on $F$ without assuming convexity.
{\renewcommand{\arraystretch}{1.5}
\begin{table}[H]
\centering
 \begin{tabularx}{\linewidth}{|>{\centering}X|>{\centering}X|>{\centering}X|}
  \hline
  \textbf{Reference} & \textbf{Assumption on $F$} & \textbf{Exponential rate $K$ of $F(x(t))-F^*$} \tabularnewline
  \hline
  Polyak \cite{polyak1964some} & $\mathcal{S}_\mu$, $C^2$ and convexity & $2\sqrt{\mu}$\tabularnewline
  \hline
  Siegel \cite{siegel2019accelerated} & $\mathcal{S}_\mu$ and convexity & $\sqrt{\mu}$ \tabularnewline
  \hline
  Aujol et al. \cite{aujol2022convergenceL} & $\mathcal{G}^2_\mu$ and convexity\\Uniqueness of the minimizer & $(2-\sqrt{2})\sqrt{\mu}$ \tabularnewline
  \hline\hline
  Polyak, Shcherbakov \cite{polyakshch} & $C^2$, \L{}ojasiewicz with $\theta=\frac{1}{2}$ and constant $c_\ell$, $L$-Lipschitz gradient & $2\frac{\mu^\frac{3}{2}}{(\sqrt{2}+1)\mu+L}$ \\with $\mu=\frac{1}{2c_\ell^2}$\tabularnewline\hline
  Apidopoulos et al. \cite{apidopoulos2022convergence} & \L{}ojasiewicz with $\theta=\frac{1}{2}$ and constant $c_\ell$, $L$-Lipschitz gradient & $2\left(\sqrt{\frac{L}{\mu}}-\sqrt{\frac{L-\mu}{\mu}}\right)\sqrt{\mu}$ \\with $\mu=\frac{1}{2c_\ell^2}$\tabularnewline\hline
  Kassing, Weissmann \cite{kassing2024polyak} & \L{}ojasiewicz with $\theta=\frac{1}{2}$ and constant $c_\ell$, $L$-Lipschitz gradient, four times differentiable & $2\sqrt{\mu}$ \\with $\mu=\frac{1}{2c_\ell^2}$\tabularnewline\hline
\end{tabularx}
\caption{Convergence rate of $F(x(t))-F^*$ where $x$ is solution of \eqref{eq:HB_ODE} for the best choice of $\alpha$ (which depends on the assumptions satisfied by $F$). The constant $K$ is defined such that $F(x(t))-F^*=\mathcal{O}\left(e^{-Kt}\right)$.}
\label{tab:hb_nu}
\end{table}}
The goal of the next part is to show that under convexity and quadratic growth conditions, a faster exponential rate, independent of $L$, can be achieved for the solution of the Heavy Ball dynamical system \eqref{eq:HB_ODE} without assuming the uniqueness of the minimizer. 

\subsection{Fast exponential decay under quadratic growth conditions}

We consider the Heavy Ball Friction (HBF) system defined as follows:
\begin{equation}
\forall t\geqslant t_0,\quad\ddot{x}(t)+\alpha\dot{x}(t)+\nabla F(x(t))=0,
\label{eq:HB_ODE}
\tag{HBF}
\end{equation}
where $t_0>0$, $\alpha>0$ and $F:\R^N\rightarrow\R$ is a convex differentiable function satisfying some quadratic growth condition. 
Generalizing recent works \cite{siegel2019accelerated,aujol2022convergenceL,aujol2022convergenceQC} and making assumptions about the regularity of the boundary of the set of minimizers, we prove that the trajectories of \eqref{eq:HB_ODE} can achieve a fast exponential decay:
\begin{theoreme}\label{thm:continu}
Let $F$ be a convex differentiable function having a non empty set of minimizers $X^*$. Suppose that $X^*$ has a $C^2$ bound or that it is a polyhedral set.
Let $x$ be a solution of \eqref{eq:HB_ODE} for some $t_0\geqslant0$ and $\alpha>0$. If $F$ satisfies the assumption $\mathcal{G}^2_\mu$ for some $\mu>0$ and $\alpha = \left(2-\frac{\sqrt{2}}{2}\right)\sqrt{\mu}$, then
\begin{equation}
\forall t\geqslant t_0,~F(x(t))-F^*\leqslant \left(\frac{11}{2}-2\sqrt{2}\right)M_0e^{-(2-\sqrt{2})\sqrt{\mu}(t-t_0)},
\end{equation}
where $M_0=F(x(t_0))-F^*+\frac{1}{2}\|\dot{x}(t_0)\|^2$. Moreover,
\begin{equation}
\|\dot{x}(t)\|=\mathcal{O}\left(e^{-\left(1-\frac{\sqrt{2}}{2}\right)\sqrt{\mu}t}\right).
\end{equation}
\end{theoreme}

We give elements of proof in the following section and a demonstration of the second claim is given in Section \ref{sec:proof_traj_hbf}.

\begin{proposition}\label{PropContinue}
Let $F$ be a convex differentiable function having a non empty set of minimizers $X^*$. Suppose that $X^*$ has a $C^2$ bound or that it is a polyhedral set. Assume that $F$ is a $\mu$-quasi strongly convex function, i.e there exists $\mu>0$ such that: 
\begin{equation*}
\forall x\in\R^N,~\langle \nabla F(x),x-x^*\rangle\geqslant F(x)-F^*+\frac{\mu}{2}\|x-x^*\|^2,
\end{equation*}
where $x^*$ denotes the projection of $x$ onto $X^*$. Let $x$ be a solution of \eqref{eq:HB_ODE} for some $t_0\geqslant0$ and $\alpha>0$. Then if $\alpha= \frac{3}{\sqrt{2}}\sqrt{\mu}$:
\begin{equation}
\forall t\geqslant t_0,~F(x(t))-F^*\leqslant 39M_0e^{-\sqrt{2\mu}(t-t_0)},
\end{equation}
where $M_0=F(x(t_0))-F^*+\frac{1}{2}\|\dot{x}(t_0)\|^2$.

\end{proposition}

\begin{remarque}
    Theorem \ref{thm:continu} and Proposition \ref{PropContinue} are based on the assumption that the set of minimizers $X^*$ has a $C^2$ bound or is a polyhedral set. More generally, the corresponding statements hold if the set $X^*$ is second order regular by the definition of Bonnans et al. \cite{bonnans1998sensitivity}, which is a weaker assumption. Given the technical nature of this hypothesis, the results are given for special cases. We refer the careful reader to the above reference for more details.
\end{remarque}

\begin{remarque}
    The fact that a regularity assumption on the set of minimizers $X^*$ is needed to obtain these results is a curiosity, since no such hypothesis is required in the discrete case, i.e. for Theorems \ref{thm:V-FISTA} and \ref{thm:vfista_v2}. It is directly related to the time-continuous nature of the trajectory $x$.
\end{remarque}

\subsubsection{Comparisons and comments}
The first study of \eqref{eq:HB_ODE} has been proposed by Polyak \cite{polyak1964some}. In this seminal work, Polyak consider a $C^2$ $\mu$-strongly convex functions. Polyak proved that for such functions the solution $x$ of \eqref{eq:HB_ODE} satisfies $F(x(t))-F(x^*)=O(e^{-2\sqrt{\mu}t})$ for a suitable choice of the friction parameter $\alpha$. If the function $F$ is $C^1$ and $\mu$-strongly convex the convergence rate is weaker, see for example \cite{siegel2019accelerated,aujol2022convergenceQC}. If $F$ is $C^1$, satisfies a quadratic growth condition and has a unique minimizer, which is a weaker assumption than strong convexity, Aujol et al. \cite{aujol2022convergenceQC} proved that the solution of \eqref{eq:HB_ODE} satisfies 
$F(x(t))-F^*=O(e^{-(2-\sqrt{2})\sqrt{\mu}t})$ for another choice of the friction parameter $\alpha$, which is slighlty slower that the rate achieved by Polyak. All the above works use the fact that the function $F$ to minimize has a unique minimizer.
Indeed if $F$ is $C^1$, convex and satisfies a quadratic growth and has several minimizers, there were no results ensuring that the solution of \eqref{eq:HB_ODE} satisfies 
$F(x_n)-F^*=0(e^{-C\sqrt{\mu}t})$ for any $C>0$. As far as we know Theorem \ref{thm:continu} is the first one ensuring such decay rate on this set of convex functions. This fast decay allows to prove the Theorem \ref{thm:vfista_v2} ensuring a fast decay of an inertial algorithm on the same set of convex functions.

Several other articles provide interesting results decay rate of the solution of the Heavy Ball ODE \eqref{eq:HB_ODE}. In \cite{polyakshch,apidopoulos2022convergence,begout2015damped,kassing2024polyak} authors provide general analysis considering \L ojasiewicz properties. In these three articles, some results on the trajectory $x(t)$ or the error $F(x(t))-F^*$ are given. It is not simple to perform a fair comparison between these results and Theorem \ref{thm:continu} because our analysis relies on the convexity and the global analysis of these works do not use this assumption. Nevertheless, in \cite{begout2015damped} and \cite{apidopoulos2022convergence} provides some decay bounds when the convexity assumption is added. More precisely, in \cite{begout2015damped}, Corollary 5.5 ensures that if $F$ is convex, $C^2$ and satisfies a quadratic growth condition with parameter $\mu$ then the trajectory $x(t)$ converges to a minimizer $x^*$of $F$, the length of the trajectory is finite and $\norm{x(t)-x^*}^2=O(e^{-\mu t})$.  Indeed, for such functions 
$d(x,X^*)^2\leqslant \frac{2}{\mu}(F(x(t))-F^*)$ and 
the Theorem \ref{thm:continu} ensures a better decay rate of the trajectory to the set of minimizers. 

The work of Apidoupoulos et al. \cite{apidopoulos2022convergence} deepens the one of Polyak et al. \cite{polyakshch} providing explicit decay of $F(x(t))-F^*$ under similar hypothesis i.e 
\L ojasiewicz properties, $C^2$ assumptions and a uniform bound on the Hessian of $F$ in the neighborhood of the set of minimizers. 
That is why we compare our results to those in \cite{apidopoulos2022convergence}, but the same conclusions hold with \cite{polyakshch}.
The bounds provided by the authors depend on a uniform bound $L$ on the Hessian of $F$ which is not the case for Theorem \ref{thm:continu} whose bound is better than Theorem 2 of \cite{apidopoulos2022convergence} when $\frac{L}{\mu}>3$. It turns out that the analysis of Apidopoulos et al has been developed in a non convex setting and in this setting, the use of this bound on the Hessian seems the only known way to get bounds on $F(x_n)-F^*$. The convexity seems to be a price to pay to get bounds independent of this Lipschitz constant $L$.


This analysis of the Heavy Ball dynamical provides a guideline for the analysis of the analysis of the corresponding discrete scheme. 

\begin{remarque}
Even if the convexity of $F$ seems to be a key hypothesis to reach such decay rate, the careful reader may notice that Theorem \ref{thm:continu} actually holds for star convex functions i.e functions satisfying: 
\begin{equation*}
\forall x\in\R^n,~\forall x^*\in X^*,~F(x)-F^*\leqslant \ps{x-x^*}{\nabla F(x)}.
\end{equation*}
where $X^*$ denotes the set of minimizers of $F$.
\end{remarque}

\subsubsection{Elements of proof}
\label{sec:discussion_NU}
The results obtained in this paper rely on a Lyapunov approach. Let us recall that when $F$ has a unique minimizer i.e $X^*=\{x^*\}$, then a classical Lyapunov choice for \eqref{eq:HB_ODE} is:
\begin{equation}
\mathcal{E}(t)=F(x(t))-F^*+\frac{1}{2}\|\lambda(x(t)-x^*)+\dot{x}(t)\|^2+\frac{\xi}{2}\|x(t)-x^*\|^2,
\end{equation}
for some well-chosen parameters $\lambda>0$ and $\xi\in\R$. Our approach to extend that type of analysis without the uniqueness assumption is to adapt the Lyapunov energy to our relaxed setting. Let $F$ have a non empty set of minimizers $X^*$ which is not reduced to one point. Let
\begin{equation}
\mathcal{E}^*(t)=F(x(t))-F^*+\frac{1}{2}\|\lambda(x(t)-x^*(t))+\dot{x}(t)\|^2+\frac{\xi}{2}\|x(t)-x^*\|^2,
\end{equation}
where for all $t\geqslant t_0$, $x^*(t)$ denotes the projection of $x(t)$ onto $X^*$, i.e.
$$x^*(t)=\arg \inf\limits_{x^*\in X^*}\|x(t)-x^*\|^2 :=P_{X^*}(x(t)).$$
This slight modification leads to a question when attempting to differentiate the Lyapunov energy: is $t\mapsto x^*(t)$ differentiable ? 

The smoothness of $t\mapsto x^*(t)$ is related to the smoothness of $P_{X^*}$. In fact, if $P_{X^*}$ is directionally differentiable then $t\mapsto x^*(t)$ is right-differentiable (and left-differentiable) and its right-hand derivative is equal to $P^\prime_{X^*}(x(t),\dot{x}(t))$. 

In \cite[Theorem~7.2]{bonnans1998sensitivity}, Bonnans et al. prove that if a closed convex set $\mathcal{S}\subset\mathcal{X}$ is second order regular at $P_\mathcal{S}(x)$ for some $x\in\mathcal{X}$, then $P_\mathcal{S}$ is directionally differentiable at $x$. We refer the reader to \cite{bonnans1998sensitivity,shapiro2016differentiability} to have a complete understanding of the notion of second order regularity. Note that this geometry assumption is satisfied by sets having a $C^2$ bound \cite{hiriart1982points} and polyhedral sets \cite{shapiro2016differentiability}. 

\begin{figure}[h]
\begin{center}
\begin{tikzpicture}
	\draw [fill=gray!20,domain=-2.5:2.5] plot (-\x*\x/3,\x) ;
	\draw (-1.2,-1.5) node[left]{$X^*$} ;
	\draw (2/3,5/3) node[right] {$x(t)$} node{$\bullet$} ;
	\draw (-1/3,1) node[below left] {$x^*(t)$} node{$\bullet$};
	\draw [<-,red,very thick] (1/6,2) -- (2/3,5/3) ;
	\draw [red] (1/6,2) node[above]{\small$\dot{x}(t)$\normalsize};
	\draw [<-,red,very thick] (-2/3,1.5) -- (-1/3,1) ;
	\draw [red] (-2/3,1.5) node[below left]{\small$\dot{x^*}(t)$\normalsize};
	\draw [->,blue,very thick] (-1/3,1) -- (2/3,5/3);
	\draw [blue] (1,0.9) node[]{\small$x(t)-x^*(t)$\normalsize};
\end{tikzpicture}
\begin{tikzpicture}
	\draw [fill=gray!20,domain=-2.5:2.5] plot (\x,-{abs(\x)}) ;
	\draw (0.6,-1.9) node[right]{$X^*$} ;
	\draw (-0.3,0.8) node[above right] {$x(t)$} node{$\bullet$} ;
	\draw (0,0) node[right] {$x^*(t)$} node{$\bullet$};
	\draw [<-,red,very thick] (-0.6,0.4) -- (-0.3,0.8) ;
	\draw [red] (-0.6,0.4) node[above left]{\small$\dot{x}(t)$\normalsize};
	\draw [red] (0,0) node[below left]{\small$\dot{x^*}(t)=0$\normalsize};
	\draw [->,blue,very thick] (0,0) -- (-0.3,0.8);
	\draw [blue] (0,0.5) node[right]{\small$x(t)-x^*(t)$\normalsize};
\end{tikzpicture}
\end{center}
\caption{Behaviour of $\dot{x^*}$ for a set of minimizers having a $C^2$ bound (on the left) and a polyhedral set of minimizers (on the right).}
\label{fig:minimizers}
\end{figure}
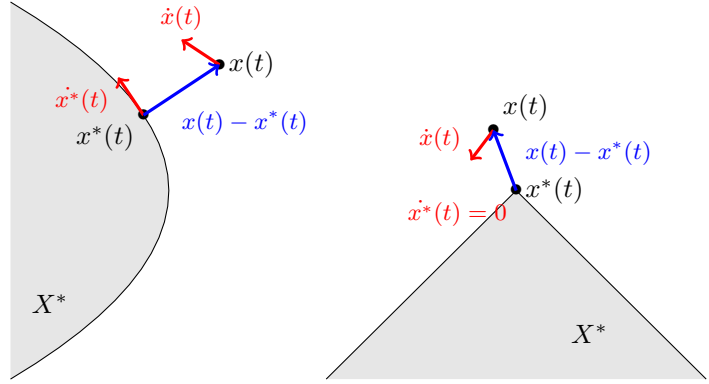

Throughout the rest of this section we assume that the set of minimizers is second order regular. Consequently, $t\mapsto x^*(t)$ is right-differentiable as well as $\mathcal{E}$. For the sake of simplicity, let $\dot{x^*}$ and $\dot{\mathcal{E}}$ denote the corresponding right-hand derivatives. We can write that:
\begin{equation}
\dot{\mathcal{E}}^*(t)=D(t)-(\lambda^2+\xi)\langle x(t)-x^*(t),\dot{x^*}(t)\rangle-\lambda\langle\dot{x}(t),\dot{x^*}(t)\rangle,
\end{equation}
where \small$$D(t)=\langle\nabla F(x(t)),\dot{x}(t)\rangle+\langle\lambda(x(t)-x^*(t))+\dot{x}(t),\lambda\dot{x}(t)+\ddot{x}(t)\rangle+\xi\langle x(t)-x^*(t),\dot{x}(t)\rangle.$$\normalsize
Observe that $D$ is exactly equal to respectively $\dot{\mathcal{E}}$ if $F$ has a unique minimizer $x^*$. The objective is then to control the additional terms $\langle x(t)-x^*(t),\dot{x^*}(t)\rangle$ and $\langle \dot{x}(t),\dot{x^*}(t)\rangle$. We introduce Figure \ref{fig:minimizers} to give an intuition of the behaviour of these terms.

We can first prove that $\langle \dot{x}(t),\dot{x^*}(t)\rangle$ is positive by using the expression $\dot{x^*}(t)=\lim\limits_{h\rightarrow0}\frac{x^*(t+h)-x^*(t)}{h}$ and the property of the projection onto a convex set. Indeed, as $X^*$ is a closed convex set, for any $x\in\R^N$ and $u\in X^*$:
\begin{equation*}
\langle x-P_{X^*}(x),u-P_{X^*}(x)\rangle\leqslant0.
\end{equation*}
Thus, for any $h>0$ we have:\small
\begin{align*}
\langle x(t+h)-x(t),x^*(t+h)-x^*(t)\rangle&=\langle x(t+h)-x^*(t+h),x^*(t+h)-x^*(t)\rangle\\&~~~+\|x^*(t+h)-x^*(t)\|^2\\&~~~+\langle x(t)-x^*(t),x^*(t)-x^*(t+h)\rangle\geqslant0.
\end{align*}\normalsize
By considering $h$ tending towards $0$ we can deduce that $\langle \dot{x}(t),\dot{x^*}(t)\rangle\geqslant0$.

In \cite[Theorem~7.2]{bonnans1998sensitivity} the authors give an expression of the directional derivative $P'_\mathcal{S}(x,d)$ for a closed convex set $\mathcal{S}\subset\mathcal{X}$ being second order regular at $P_\mathcal{S}(x)$ for some $x\in\mathcal{X}$. This directional derivative satisfies:
\begin{equation*}
\langle x-P_\mathcal{S}(x),P'_\mathcal{S}(x,d)\rangle=0.
\end{equation*}
Considering the assumptions made on $X^*$ we can deduce that $\langle x(t)-x^*(t),\dot{x^*}(t)\rangle=0$ for all $t\geqslant t_0$. 

These results ensure that for any choice of parameter $\lambda>0$ and $\xi\in\R$,$
\dot{\mathcal{E}}^*(t)\leqslant D(t)$. From this point, the proofs of the convergence results stated in Theorem \ref{thm:continu} and Proposition \ref{PropContinue} follow the original proofs, taking $D$ instead of $\dot{\mathcal{E}}$ and by applying the following lemma. A proof is given in Section \ref{sec:proof_right_diff}.

\begin{lemme}
Let $\phi:\R\rightarrow\R$ be a continuous function which is right-differentiable. Assume that 
\begin{equation}
\forall t\geqslant t_0,~\phi_+(t)\leqslant\psi(t),
\label{eq:ineq_phi+}
\end{equation}
where $\phi_+(t)=\lim\limits_{h\rightarrow0,~h>0}\frac{\phi(t+h)-\phi(t)}{h}$ denotes the right derivative of $\phi$ at $t$. Then,
\begin{equation}
\forall t\geqslant t_0,~\phi(t)\leqslant\phi(t_0)+\int_{t_0}^t\psi(u)du.
\label{eq:maj_semidif}
\end{equation}
\label{lem:right-diff}
\end{lemme}

\section{Proofs of Theorems~\ref{thm:V-FISTA} and \ref{thm:vfista_v2} }\label{sec:proofs}
\label{sec:proof_notations}

The proofs of Theorems \ref{thm:V-FISTA} and \ref{thm:vfista_v2} are built around the following discrete Lyapunov energies
\begin{equation}
\mathcal{E}_n=\frac{2}{L}(F(x_n)-F^*)+\left\|\lambda(x_{n-1}-x_{n-1}^*)+x_n-x_{n-1}\right\|^2,
\end{equation}
and
\begin{equation}
\mathcal{E}_n=\frac{2}{L}(F(x_n)-F^*)+\alpha\left\|\lambda(x_{n}-x_{n}^*)+x_n-x_{n-1}\right\|^2+\lambda(1-\alpha)^2\|x_n-x_n^*\|^2,
\end{equation}
where $\lambda>0$ and $x_n^*$ denotes the projection of $x_n$ onto $X^*$ for any $n\in\N$. For the sake of clarity, we use the following notations:
\begin{equation}
\begin{gathered}
w_n=\frac{2}{L}(F(x_n)-F^*)),~h_n=\|x_n-x^*_n\|^2,\\
\delta_n=\|x_n-x_{n-1}\|^2,~\gamma_n^*=\|x_n^*-x_{n-1}^*\|^2.
\end{gathered}\label{eq:notations}
\end{equation}
The following lemma is necessary in order to handle the terms related to non uniqueness of the minimizers. We give a proof in Section \ref{sec:proof_lemma_tech1}.
\begin{lemme}
For all $n\in\N^*$, the following equalities hold:
\begin{enumerate}
	\item $\langle x_n-x_n^*,x_n-x_{n-1}\rangle=\frac{1}{2}(h_n-h_{n-1}+\delta_n-\gamma_n^*)+\langle x_{n-1}-x_{n-1}^*,x_n^*-x_{n-1}^*\rangle.$
	\item $\langle x_{n-1}-x_{n-1}^*,x_n-x_{n-1}\rangle=\frac{1}{2}(h_n-h_{n-1}-\delta_n+\gamma_n^*)+\langle x_{n}-x_{n}^*,x_n^*-x_{n-1}^*\rangle.$
\end{enumerate}
\label{lem:tech1}
\end{lemme}
Moreover, we introduce a lemma which encodes the fact that the sequence $(x_n)_{n\in\N}$ is provided by \eqref{eq:V-FISTA}. The proof is based on the descent lemma proved in \cite{chambolle2015convergence} and it can be found in Appendix \ref{sec:proof_lemma_tech2}.
\begin{lemme}
Let $(x_n)_{n\in\N}$ be the sequence provided by \eqref{eq:V-FISTA} with $s=\frac{1}{L}$. Then, for any $n\in\N^*$,
\begin{eqnarray*}
w_{n+1}-w_n&\leqslant& \alpha^2\delta_n-\delta_{n+1},\\
w_{n+1}&\leqslant& (1+\alpha)h_n+(\alpha^2+\alpha)\delta_n-\alpha h_{n-1}-h_{n+1}-\gamma_{n+1}^*-\alpha\gamma_n^*\vspace{.2cm}\\
&&+2\alpha\langle x_{n-1}-x_{n-1}^*,x_n^*-x_{n-1}^*\rangle-2\langle x_{n+1}-x_{n+1}^*,x_{n+1}^*-x_n^*\rangle.\nonumber
\end{eqnarray*}
\label{lem:tech2}
\end{lemme}
We would like to point out that several controls are deduced from the properties of the projection onto a convex. Indeed, if $C$ is a closed convex set such that $C\subset \R^N$, then for any $x\in \R^N$ and $y\in C$,
$$\langle x-p,y-p\rangle\leqslant0,$$
where $p$ denotes the projection of $x$ onto $C$. This property directly guarantees inequalities such as
$$\langle x_{n}-x_{n}^*,x_n^*-x_{n-1}^*\rangle\geqslant0,$$
or
$$\langle x_{n-1}-x_{n-1}^*,x_n^*-x_{n-1}^*\rangle\leqslant0.$$
\subsection{Proof of Theorem \ref{thm:V-FISTA}}
\label{sec:proof_V-FISTA}
Let $(x_n)_{n\in\N}$ be the sequence provided by \eqref{eq:V-FISTA} for some $\alpha>0$ to be defined. We define the following discrete Lyapunov energy:
\begin{equation}
\mathcal{E}_n=\frac{2}{L}(F(x_n)-F^*)+\|x_n-x_{n-1}+\lambda(x_{n-1}-x_{n-1}^*)\|^2,
\label{eq:Lyap_VFISTA}
\end{equation}
where $\lambda>0$ and $x_n^*$ denotes the projection of $x_n$ on $X^*$ for any $n\in\N$. By setting $\lambda=\sqrt{\kappa}$ and considering that $F$ has a unique minimizer, we recover the energy considered by Beck in \cite{beck2017first}. 

The aim of this proof is to find $\tau>0$ as large as possible such that for a well-chosen set of parameters,
\begin{equation}
\mathcal{E}_{n+1}-(1-\tau\sqrt{\kappa})\mathcal{E}_n\leqslant0.
\label{eq:wanted_ineq}
\end{equation}

The proof is divided into three parts. We first use the lemmas introduced in the introduction of Section \ref{sec:proof_notations} and the properties of the projection onto a convex to handle the terms related to the non uniqueness of the minimizers. Then, we give a set of parameters which leads to the wanted inequality \eqref{eq:wanted_ineq} by using the geometry assumption satisfied by $F$. The convergence of the trajectories is obtained in the last section using the previous results and elementary computations.

\subsubsection{Preliminary work}

We recall that we use the notations defined in \eqref{eq:notations}. By rewriting \eqref{eq:Lyap_VFISTA} and using the second claim of Lemma \ref{lem:tech1} we have:
\begin{equation}
\mathcal{E}_n=w_n+(1-\lambda)\delta_n+\lambda (h_n-h_{n-1})+\lambda^2h_{n-1}+\lambda\gamma_n^*+2\lambda\langle x_n-x_n^*,x_n^*-x_{n-1}^*\rangle.
\label{eq:Lyap_rewr}
\end{equation}
Lemma \ref{lem:tech2} ensures that if $\lambda\leqslant1$:
\begin{equation}
\begin{aligned}
w_{n+1}-(1-\lambda)w_n&\leqslant\alpha(\alpha+\lambda)\delta_n-(1-\lambda)\delta_{n+1}+\alpha\lambda(h_n- h_{n-1})+\lambda(h_n- h_{n+1})\\&-\lambda\gamma_{n+1}^*-\alpha\lambda\gamma_n^*+2\alpha\lambda\langle x_{n-1}-x_{n-1}^*,x_n^*-x_{n-1}^*\rangle\\&-2\lambda\langle x_{n+1}-x_{n+1}^*,x_{n+1}^*-x_n^*\rangle.
\end{aligned}
\end{equation}
This inequality combined with \eqref{eq:Lyap_rewr} ensures that:
\begin{equation}
\mathcal{E}_{n+1}-(1-\lambda)\mathcal{E}_n\leqslant a_1\delta_n+a_2(h_n-h_{n-1})+a_3 h_n+\mathcal{X}_n^*,
\end{equation}
where:
$$a_1=\alpha(\alpha+\lambda)-(1-\lambda)^2,~a_2=\alpha\lambda-\lambda(1-\lambda)+(1-\lambda)\lambda^2,~a_3=\lambda^3,$$
and $\mathcal{X}_n^*$ is defined by
\begin{align*}\mathcal{X}_n^*=&-\lambda(1-\lambda+\alpha)\gamma_n^*+2\alpha\lambda\langle x_{n-1}-x_{n-1}^*,x_n^*-x_{n-1}^*\rangle\\&-2\lambda(1-\lambda)\langle x_n-x_n^*,x_n^*-x_{n-1}^*\rangle.\end{align*}
Due to the properties of the projection onto a convex set we have that for all $n\in\N$:
\begin{equation*}
\left\{
\begin{gathered}
\langle x_{n-1}-x_{n-1}^*,x_n^*-x_{n-1}^*\rangle\leqslant0,\\
\langle x_n-x_n^*,x_n^*-x_{n-1}^*\rangle\geqslant0.
\end{gathered}\right.
\end{equation*}
and since $\gamma_n^*\geqslant0$, we can conclude that $\mathcal{X}_n^*\leqslant0$ and consequently that:
\begin{equation}
\mathcal{E}_{n+1}-(1-\lambda)\mathcal{E}_n\leqslant a_1\delta_n+a_2(h_n-h_{n-1})+a_3 h_n.
\end{equation}

\subsubsection{Getting the convergence rate}

Recall that we want to find $\tau>0$ such that: $\mathcal{E}_{n+1}-(1-\tau\sqrt{\kappa})\mathcal{E}_n\leqslant0$. We choose the following set of parameters: $$\tau=\frac{2}{3\sqrt{3}},\quad\lambda=\frac{1}{\sqrt{3}}\sqrt{\kappa},\quad\alpha=1-\frac{5}{3\sqrt{3}}\sqrt{\kappa}=1-\frac{5}{3}\lambda.$$ Then we get that:
\begin{equation}
\begin{aligned}
\mathcal{E}_{n+1}-(1-\tau\sqrt{\kappa})\mathcal{E}_n&=\mathcal{E}_{n+1}-(1-\lambda)\mathcal{E}_n+\left(\frac{2}{3\sqrt{3}}-\frac{1}{\sqrt{3}}\right)\sqrt{\kappa}\mathcal{E}_n\\
&\leqslant a_1\delta_n+a_2(h_n-h_{n-1})+a_3 h_n-\frac{1}{3\sqrt{3}}\sqrt{\kappa}\mathcal{E}_n,
\end{aligned}
\end{equation}
where for this parameter choice we have:
\begin{eqnarray*}
    a_1 &=& -\frac{\lambda}{3}\left(1-\frac{\lambda}{3}\right)=\frac{\sqrt{\kappa}}{27}(\sqrt{\kappa}-3\sqrt{3}),\quad a_3= \lambda^3=\frac{\kappa^{\frac{3}{2}}}{3\sqrt{3}},\\
    a_2&=&\lambda^2\left(\frac{1}{3}-\lambda\right)=\frac{\kappa}{3\sqrt{3}}\left(\frac{1}{\sqrt{3}}-\sqrt{\kappa}\right).
\end{eqnarray*}
Under the condition $\kappa\leqslant \frac{1}{3}$, we have that $a_1\leqslant0$ and hence,
\begin{equation}
\mathcal{E}_{n+1}-\left(1-\frac{2}{3\sqrt{3}}\sqrt{\kappa}\right)\mathcal{E}_n\leqslant\frac{\kappa}{3\sqrt{3}}\left(\frac{1}{\sqrt{3}}-\sqrt{\kappa}\right)(h_n-h_{n-1})+\frac{\kappa^{\frac{3}{2}}}{3\sqrt{3}}h_n-\frac{1}{3\sqrt{3}}\sqrt{\kappa}\mathcal{E}_n.
\end{equation}

Moreover, as the condition $\kappa\leqslant\frac{1}{3}$ ensures that $a_2\geqslant0$ we can apply the following lemma which is proved in Section \ref{sec:proof_lem_nrj}
\begin{lemme}
Let $(x_n)_{n\in\N}$ be the sequence provided by \eqref{eq:V-FISTA} and $\lambda=\frac{1}{\sqrt{3}}\sqrt{\kappa}$. Then for all $n\in\N$:
\begin{equation}
h_n-h_{n-1}\leqslant \frac{\sqrt{3}}{\sqrt{\kappa}}\left(\mathcal{E}_n-w_n\right).
\end{equation}
\label{lem:VFISTA_nrj}
\end{lemme}
Lemma \ref{lem:VFISTA_nrj} guarantees that if $\kappa\leqslant\frac{1}{3}$:
\begin{equation}
\mathcal{E}_{n+1}-(1-\frac{2}{3\sqrt{3}}\sqrt{\kappa})\mathcal{E}_n\leqslant \frac{\kappa^{\frac{3}{2}}}{3\sqrt{3}}h_n-\frac{\sqrt{\kappa}}{3}\left(\frac{1}{\sqrt{3}}-\sqrt{\kappa}\right)w_n-\frac{\kappa}{3}\mathcal{E}_n.
\end{equation}
Moreover, as $F$ satisfies $\mathcal{G}^2_\mu$ we can write that $h_n\leqslant \frac{w_n}{\kappa}$ and consequently:
\begin{equation}
\begin{aligned}
\mathcal{E}_{n+1}-(1-\frac{2}{3\sqrt{3}}\sqrt{\kappa})\mathcal{E}_n&\leqslant \left(\frac{\sqrt{\kappa}}{3\sqrt{3}}-\frac{\sqrt{\kappa}}{3}\left(\frac{1}{\sqrt{3}}-\sqrt{\kappa}\right)\right)w_n-\frac{\kappa}{3}\mathcal{E}_n\\
&\leqslant \frac{\kappa}{3}w_n-\frac{\kappa}{3}\mathcal{E}_n.
\end{aligned}
\end{equation}
Noticing that $w_n\leqslant \mathcal{E}_n$ we can conclude that:
\begin{equation}
\mathcal{E}_{n+1}-(1-\frac{2}{3\sqrt{3}}\sqrt{\kappa})\mathcal{E}_n\leqslant0.
\label{eq:last_step_vfista}
\end{equation}
Hence: 
$\mathcal{E}_n\leqslant \left(1-\frac{2}{3\sqrt{3}}\sqrt{\kappa}\right)^n\mathcal{E}_0.$
As we consider that $x_{-1}=x_0$, we have that $\mathcal{E}_0=w_0+\lambda^2h_0$. As a consequence, the geometry condition $\mathcal{G}^2_\mu$ ensures that $\mathcal{E}_0\leqslant \frac{4}{3}w_0$ and
\begin{equation}
F(x_n)-F^*\leqslant \frac{4}{3}\left(1-\frac{2}{3\sqrt{3}}\sqrt{\kappa}\right)^n(F(x_0)-F^*).
\label{eq:rate_sigma0}
\end{equation}

\subsubsection{Convergence of the trajectories}

Let $b_n=\|x_n-x_{n-1}+\lambda(x_{n-1}-x_{n-1}^*)\|^2$. Using the inequality $\|u\|^2=2\|u+v\|^2+2\|v\|^2$, we get that:
\begin{equation}
\delta_n\leqslant2b_n+\frac{2}{\lambda^2}h_{n-1}.
\end{equation}
Thus, using the definition of $\mathcal{E}$ and the geometry of $F$:
\begin{equation}
\delta_n\leqslant 2\mathcal{E}_n+\frac{2}{\lambda^2\kappa}w_{n-1}.
\end{equation}
Then, by applying inequality \eqref{eq:last_step_vfista}, we deduce that:
\begin{equation}
\delta_n\leqslant \left(2\left(1-\frac{2}{3\sqrt{3}}\sqrt{\kappa}\right)+\frac{2}{\lambda^2\kappa}\right)\mathcal{E}_{n-1},
\end{equation}
and consequently,
\begin{equation}
\delta_n\leqslant \left(2\left(1-\frac{2}{3\sqrt{3}}\sqrt{\kappa}\right)+\frac{2}{\lambda^2\kappa}\right)\left(1-\frac{2}{3\sqrt{3}}\sqrt{\kappa}\right)^{n-1}\mathcal{E}_0.
\end{equation}
Hence,
\begin{equation}
\|x_n-x_{n-1}\|=\mathcal{O}\left(e^{-\frac{1}{3\sqrt{3}}\sqrt{\kappa}n}\right).
\end{equation}

\subsection{Proof of Theorem \ref{thm:vfista_v2}}
\label{sec:proof_v-fista2}

\subsubsection{Structure of the proof}

The proof of Theorem \ref{thm:vfista_v2} is built around the approach provided in \cite{aujol2022convergenceL} in order to prove convergence rates of the trajectories of the Heavy Ball system described by:
\begin{equation}
\ddot{x}(t)+\alpha_c\dot{x}(t)+\nabla F(x(t))=0,
\tag{HBF}
\label{eq:HBF}
\end{equation}
for some $\alpha_c>0$. Indeed, the sequence generated by \eqref{eq:V-FISTA} can be seen as a discretization of \eqref{eq:HBF} (when $F$ is differentiable) and the strategy can be adapted to the discrete setting. Note that the parameter $\alpha_c$ in \eqref{eq:HBF} does not play the same role as $\alpha$ in \eqref{eq:V-FISTA}. Indeed, we have that $\alpha_c$ behaves as $\sqrt{L}(1-\alpha)$. \\
This proof relies on the analysis of the following Lyapunov energy:
\begin{equation}
\forall n\in\N,\quad \mathcal{E}_n=\frac{2}{L}(F(x_n)-F^*)+\alpha\|x_n-x_{n-1}+\lambda(x_{n}-x_{n}^*)\|^2+\lambda(1-\alpha)^2\|x_n-x_n^*\|^2.
\end{equation}
The strategy of the proof is straightforward: we aim to find a set of parameters $(\alpha,\lambda,\nu)\in\left(\R^+\right)^3$ such that for any $n\in\N$,
\begin{equation}
\mathcal{E}_{n+1}-\mathcal{E}_n+\nu\mathcal{E}_{n+1}\leqslant0.
\end{equation} 
In this way, simple calculations show that it ensures
\begin{equation}
\forall n\in\N,\quad \mathcal{E}_n\leqslant \left(1-\nu+\nu^2\right)^n\mathcal{E}_0
\end{equation}
which leads us to the conclusion.

\subsubsection{Proof}

Let $(x_n)_{n\in\N}$ be the sequence provided by \eqref{eq:V-FISTA}. Following the notations introduced in \eqref{eq:notations}, we can rewrite:
\begin{equation}
\forall n\in\N, ~\mathcal{E}_n=w_n+\lambda(\lambda\alpha+(1-\alpha)^2)h_n+\alpha\delta_n+2\alpha\lambda\langle x_n-x_n^*,x_n-x_{n-1}\rangle.
\end{equation}
Following the first claim of \ref{lem:tech1}, 
\begin{equation}
\begin{aligned}
\mathcal{E}_n=&w_n+\lambda(\lambda\alpha+(1-\alpha)^2)h_n+\lambda\alpha(h_n-h_{n-1})+(1+\lambda)\alpha\delta_n\\
&-\lambda\alpha\gamma_n^*+2\lambda\alpha\langle x_{n-1}-x_{n-1}^*,x_n^*-x_{n-1}^*\rangle.
\end{aligned}\label{eq:rewriting_lyap2}
\end{equation}
Observe that due to the properties of the projection onto a convex, for any $n\in\N$:
\begin{equation}
\mathcal{E}_n\leqslant w_n+\lambda(\lambda\alpha+(1-\alpha)^2)h_n+\lambda\alpha(h_n-h_{n-1})+(1+\lambda)\alpha\delta_n.
\label{eq:VFISTA_lyap_sb}
\end{equation}
By exploiting the expression \eqref{eq:rewriting_lyap2}, we can show the following lemma. The proof can be found in Section \ref{sec:proof_lem_diff_vf2}.
\begin{lemme}
For any $n\in\N$, we have that:
\begin{equation}
\begin{aligned}
\mathcal{E}_{n+1}-\mathcal{E}_n\leqslant &~-\lambda w_{n+1}+\lambda\alpha(\lambda+\alpha-1)(h_{n+1}-h_n)\\&+(\lambda\alpha+\alpha-1)\delta_{n+1}-\alpha(1-\alpha-\lambda\alpha)\delta_n.
\end{aligned}
\end{equation}
\label{lem:VFISTA_diff}
\end{lemme}
This inequality combined to \eqref{eq:VFISTA_lyap_sb} guarantees that for any $\nu>0$,
\begin{equation}
\begin{aligned}
\mathcal{E}_{n+1}-\mathcal{E}_n+\nu\mathcal{E}_{n+1}\leqslant &~(\nu-\lambda) w_{n+1}+\lambda\alpha(\lambda+\alpha-1+\nu)(h_{n+1}-h_n)\\&+((1+\lambda)\alpha(1+\nu)-1)\delta_{n+1}-\alpha(1-\alpha-\lambda\alpha)\delta_n\\&+\nu\lambda(\lambda\alpha+(1-\alpha)^2)h_{n+1}.
\end{aligned}
\end{equation}
We make the following choice of parameters:
\begin{equation}
\alpha=1-\omega\sqrt{\kappa},~\nu=\tau\sqrt{\kappa},~\lambda=1-\alpha-\nu=(\omega-\tau)\sqrt{\kappa}.
\end{equation}
This set of parameters ensures that the following inequality is valid:
\begin{equation}\begin{aligned}
\mathcal{E}_{n+1}-\mathcal{E}_n+\tau\sqrt{\kappa}\mathcal{E}_{n+1}\leqslant &~(2\tau-\omega)\sqrt{\kappa} w_{n+1}-(\omega\tau+(\omega-\tau)^2+\omega\tau(\omega-\tau)\sqrt{\kappa})\kappa\delta_{n+1}\\&-(1-\omega\sqrt{\kappa})(\tau+(\omega-\tau)\sqrt{\kappa})\sqrt{\kappa}\delta_n\\&+\tau(\omega-\tau)(\omega-\tau+\omega\tau\sqrt{\kappa})\kappa^\frac{3}{2} h_{n+1}.\end{aligned}
\end{equation}
Note that we consider a parameter $\alpha>0$ which implies that $1-\omega\sqrt{\kappa}>0$. Suppose in addition that $\omega>2\tau$. Then, we get that for any $n\in\N$,
\begin{equation}
\mathcal{E}_{n+1}-\mathcal{E}_n+\tau\sqrt{\kappa}\mathcal{E}_{n+1}\leqslant-(\omega-2\tau)\sqrt{\kappa}w_{n+1}+\tau(\omega-\tau)(\omega-\tau+\omega\tau\sqrt{\kappa})\kappa^\frac{3}{2} h_{n+1}.
\end{equation}
As $F$ satisfies the assumption $\mathcal{G}^2_\mu$, we can write that $\kappa h_{n+1}\leqslant w_{n+1}$ and hence:
\begin{equation}
\mathcal{E}_{n+1}-\mathcal{E}_n+\tau\sqrt{\kappa}\mathcal{E}_{n+1}\leqslant\left(\tau(\omega-\tau)(\omega-\tau+\omega\tau\sqrt{\kappa})-\omega+2\tau\right)\kappa^\frac{3}{2}h_{n+1}.
\end{equation}
Thus, if
\begin{equation}
\left(1-\omega\sqrt{\kappa}\right)\tau^3-\omega\left(2-\omega\sqrt{\kappa}\right)\tau^2+(\omega^2+2)\tau-\omega\leqslant0,\label{eq:cond_in_proof}
\end{equation}
then for any $n\in\N$,
\begin{equation}
\mathcal{E}_{n+1}-\mathcal{E}_n+\tau\sqrt{\kappa}\mathcal{E}_{n+1}\leqslant0.
\end{equation}
Note that the solutions of \eqref{eq:cond_in_proof} automatically satisfy $\omega>2\tau$. Elementary computations show that this implies
\begin{equation}
\mathcal{E}_{n}\leqslant \left(1-\tau\sqrt{\kappa}+\tau^2\kappa\right)^n\mathcal{E}_0.\label{eq:ineq_vfistav2}
\end{equation}
Note that since $x_{-1}=x_0$,
\begin{equation}
\mathcal{E}_0=w_0+\lambda\left(\lambda\alpha+(1-\alpha)^2\right)h_0=w_0+\left((\omega-\tau)^2+(\omega-\tau)\omega\tau\sqrt{\kappa}\right)\kappa h_0,
\end{equation}
and using the assumption $\mathcal{G}^2_\mu$,
\begin{equation}
\mathcal{E}_0\leqslant \left(1+(\omega-\tau)^2+(\omega-\tau)\omega\tau\sqrt{\kappa}\right)w_0.
\end{equation}
Moreover, for any $n\in\N$, $w_n\leqslant \mathcal{E}_n$. Thus, if \eqref{eq:cond_in_proof} is satisfied, then for any $n\in\N$:
\begin{equation}
F(x_n)-F^*\leqslant \left(1+(\omega-\tau)^2+(\omega-\tau)\omega\tau\sqrt{\kappa}\right)\left(1-\tau\sqrt{\kappa}+\tau^2\kappa\right)^n(F(x_0)-F^*).
\end{equation}
In addition, by applying the inequality $\|u\|^2=2\|u+v\|^2+2\|v\|^2$, we get that:
\begin{equation}
\delta_n\leqslant2b_n+\frac{2}{\lambda^2}h_{n},
\end{equation}
where $b_n=\|\lambda(x_n-x_n^*)+x_n-x_{n-1}\|^2$. The assumption $\mathcal{G}^2_\mu$ gives that
\begin{equation}
\delta_n\leqslant 2b_n+\frac{2}{\lambda^2\kappa}w_{n},
\end{equation}
and given the definition of $\mathcal{E}_n$ we obtain,
\begin{equation}
\delta_n\leqslant \left(\frac{2}{\alpha}+\frac{2}{\lambda^2\kappa}\right)\mathcal{E}_{n}.
\end{equation}
By combining the above inequality with \eqref{eq:ineq_vfistav2}, we can finally prove that if \eqref{eq:cond_in_proof} is valid, then
\begin{equation}
\|x_n-x_{n-1}\|=\mathcal{O}\left(e^{-\frac{1}{2}\tau\sqrt{\kappa}\left(1-\tau\sqrt{\kappa}\right)n}\right).
\end{equation}

\subsection{Proof of Corollary \ref{cor:vfista_slow}}\label{sec:proof_cor2}

\begin{figure}[ht]
\centering
\includegraphics[width=0.9\textwidth]{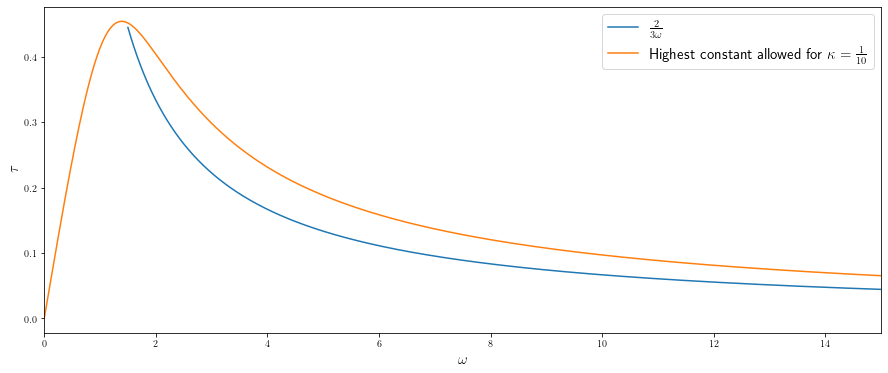}
\caption{Comparison of the highest rate $\tau$ satisfying \eqref{eq:cond_vfista_eg} for $\kappa=\frac{1}{10}$ and $\omega>0$ with the function $\omega\mapsto\frac{2}{3\omega}$.}
\label{fig:VFISTA_approx}
\end{figure}

Let $(x_n)_{n\in\N}$ be the sequence given by \eqref{eq:V-FISTA} for some $\alpha>0$ and $s=\frac{1}{L}$. According to Theorem \ref{thm:vfista_v2}, if $\alpha=1-\omega\sqrt{\kappa}$ for some $\omega\in\left(0,\frac{1}{\sqrt{\kappa}}\right)$, then \eqref{eq:thm_v2_1} and \eqref{eq:thm_v2_2} are valid for any $\tau>0$ satisfying:
\begin{equation}
\left(1-\omega\sqrt{\kappa}\right)\tau^3-\omega\left(2-\omega\sqrt{\kappa}\right)\tau^2+(\omega^2+2)\tau-\omega\leqslant0.\label{eq:cond_vfista_cor}
\end{equation}
Corollary \ref{cor:vfista_slow} relies on the following lemma.
\begin{lemme}\label{lem:maj_omega}
    Let $$P:(\tau,\omega,\kappa)\mapsto\left(1-\omega\sqrt{\kappa}\right)\tau^3-\omega\left(2-\omega\sqrt{\kappa}\right)\tau^2+(\omega^2+2)\tau-\omega.$$ If $\kappa\leqslant \frac{1}{10}$, then for any $\omega\geqslant\frac{3}{2}$, $P\left(\frac{2}{3\omega},\omega,\kappa\right)<0$. 
\end{lemme}
\noindent\textbf{Proof.}
Given the expression of $P$, simple computations give that for any $\kappa\in(0,1)$ and $\omega>0$,
\begin{equation}
    P\left(\frac{2}{3\omega},\omega,\kappa\right)=\left(1-\omega\sqrt{\kappa}\right)\frac{8-12\omega^2}{27\omega^3}+\frac{8-3\omega^2}{9\omega}.
\end{equation}
We define the function $\Phi$ as follows
\begin{equation*}
    \Phi(\omega;\kappa)=-9\omega^4+12\omega^3\sqrt{\kappa}+12\omega^2-8\omega\sqrt{\kappa}+8,
\end{equation*}
such that $\frac{\Phi(\omega;\kappa)}{27\omega^3}=P\left(\frac{2}{3\omega},\omega,\kappa\right)$. We get that:
\begin{equation}
    \frac{\partial\Phi}{\partial\omega}(\omega;\kappa)=-36\omega^2(\omega-\sqrt{\kappa})+24\omega-8\sqrt{\kappa}.
\end{equation}
Consequently, if $\omega\geqslant\frac{3}{2}$ and $\kappa\leqslant \frac{1}{10}$, we have that $\omega-\sqrt{\kappa}>1$ and
\begin{equation*}
    \frac{\partial\Phi}{\partial\omega}(\omega;\kappa)\leqslant-36\omega^2+24\omega<0.
\end{equation*}
Since $\Phi\left(\frac{3}{2};\kappa\right)=-\frac{169}{16}+\frac{57}{2}\sqrt{\kappa}$ which is strictly negative if $\kappa\leqslant \frac{1}{10}$, we can deduce that $\Phi(\omega;\kappa)<0$ for any $\omega\geqslant \frac{3}{2}$. As $\frac{\Phi(\omega;\kappa)}{27\omega^3}=P\left(\frac{2}{3\omega},\omega,\kappa\right)$, this lemma is proved.\\
In fact, $\omega\mapsto \frac{2}{3\omega}$ is a lower bound of the highest value of $\tau$ satisfying \eqref{eq:cond_vfista_eg} for some $\omega\geqslant \frac{3}{2}$ and $\kappa=\frac{1}{10}$ as illustrated in Figure \ref{fig:VFISTA_approx}.

According to Lemma \ref{lem:maj_omega}, if $\kappa\leqslant \frac{1}{10}$ and $\omega\in\left(\frac{3}{2},\frac{1}{\sqrt{\kappa}}\right)$, then \eqref{eq:thm_v2_1} and \eqref{eq:thm_v2_2} are valid with $\tau=\frac{2}{3\omega}$. Moreover, if we consider $\alpha=1-\theta$ for some $\theta\in\left[\frac{3}{2}\sqrt{\kappa},1\right)$, then \eqref{eq:thm_v2_1} and \eqref{eq:thm_v2_2} are satisfied with $\tau=\frac{2}{3\theta}\sqrt{\kappa}$. This leads to the conclusion of Corollary \ref{cor:vfista_slow}.

\appendix
\section{Technical proofs}

\subsection{Proof of Theorem \ref{thm:continu}}
\label{sec:proof_traj_hbf}
Our analysis follows that introduced in \cite{aujol2022convergenceL}. We set $\alpha=\left(2-\frac{\sqrt{2}}{2}\right)\sqrt{\mu}$ and we consider the following Lyapunov energy:
\begin{equation}
\mathcal{E}(t)=F(x(t))-F^*+\frac{1}{2}\|\lambda(x(t)-x^*(t))+\dot{x}(t)\|^2+\xi\|x(t)-x^*(t)\|^2,
\end{equation}
with $\lambda=\sqrt{\mu}$ and $\xi=-\left(1-\frac{\sqrt{2}}{2}\right)\mu$.

Following the discussion of Section \ref{sec:discussion_NU}, the assumptions of Theorem \ref{thm:continu} ensure that $\mathcal{E}$ is right-differentiable and for all $t\geqslant t_0$,\small
\begin{equation}
    \dot{\mathcal{E}}(t)\leqslant -\lambda \langle \nabla F(x(t)),x(t)-x^*(t)\rangle+(\lambda-\alpha)\|\dot{x}(t)\|^2+(\xi+\lambda(\lambda-\alpha))\langle x(t)-x^*(t),\dot{x}(t)\rangle,
\end{equation}\normalsize
where $\dot{\mathcal{E}}$ denotes the right derivative of $\mathcal{E}$. By using the convexity of $F$ and replacing the parameters by their value,\small
\begin{equation}
    \dot{\mathcal{E}}(t)\leqslant -\sqrt{\mu}(F(x(t))-F^*)-\left(1-\frac{\sqrt{2}}{2}\right)\sqrt{\mu}\|\dot{x}(t)\|^2-(2-\sqrt{2})\mu\langle x(t)-x^*(t),\dot{x}(t)\rangle.
\end{equation}\normalsize
Let us define $\delta=(2-\sqrt{2})\sqrt{\mu}$. The above inequality guarantees that for all $t\geqslant t_0$:\small
\begin{equation}
    \dot{\mathcal{E}}(t)+\delta\mathcal{E}(t)\leqslant \left(1-\sqrt{2}\right)\sqrt{\mu}\left(F(x(t))-F^*\right)+\frac{\sqrt{2}-1}{2}\mu^\frac{3}{2}\|x(t)-x^*(t)\|^2.
\end{equation}\normalsize
As $F$ satisfies $\mathcal{G}^2_\mu$ and $\|x(t)-x^*(t)\|=d(x(t),X^*)$ for all $t\geqslant t_0$, we obtain that for all $t\geqslant t_0$:
\begin{equation}
    \dot{\mathcal{E}}(t)+\delta\mathcal{E}(t)\leqslant \left(\left(1-\sqrt{2}\right)\sqrt{\mu}\frac{\mu}{2}+\frac{\sqrt{2}-1}{2}\mu^\frac{3}{2}\right)\|x(t)-x^*(t)\|^2\leqslant0.
\end{equation}
We refer the reader to the proof of \cite[Theorem~1]{aujol2022convergenceL} for further developments on each of the above steps and a discussion on the value of the parameters. Lemma \ref{lem:right-diff} then guarantees that, for all $t\geqslant t_0$,
\begin{equation}
\mathcal{E}(t)\leqslant \mathcal{E}(t_0)e^{-(2-\sqrt{2})\sqrt{\mu}(t-t_0)}.\label{eq:control_nrj_hbf}
\end{equation}
Since $F$ satisfies $\mathcal{G}^2_\mu$, elementary computations show that:
\begin{equation}
F(x(t))-F^*+\xi\|x(t)-x^*(t)\|^2\geqslant(\sqrt{2}-1)(F(x(t))-F^*),
\end{equation}
and consequently,
\begin{equation}
\forall t\geqslant t_0,~\mathcal{E}(t)\geqslant \left(\sqrt{2}-1\right)(F(x(t))-F^*)+\frac{1}{2}\|\lambda(x(t)-x^*(t))+\dot{x}(t)\|^2.
\end{equation}
This inequality implies that for all $t\geqslant t_0$:
\begin{equation}
F(x(t))-F^*\leqslant \frac{1}{\sqrt{2}-1}\mathcal{E}(t),\label{eq:control_F_cont}
\end{equation}
and
\begin{equation}
\|\lambda(x(t)-x^*(t))+\dot{x}(t)\|^2\leqslant 2\mathcal{E}(t).
\end{equation}
The first statement of Theorem \ref{thm:continu} can be demonstrated by combining \eqref{eq:control_nrj_hbf} and \eqref{eq:control_F_cont} and rewriting $\mathcal{E}(t_0)$ (see \cite[Section~6.1]{aujol2022convergenceL} for further details). We prove the second result as follows.\\
Using inequality $\|u\|^2\leqslant 2\|u+v\|^2+2\|v\|^2$, we get that:
\begin{equation}
\begin{aligned}
\|\dot{x}(t)\|^2&\leqslant2\|\lambda(x(t)-x^*(t))+\dot{x}(t)\|^2+2\mu\|x(t)-x^*(t)\|^2\\&\leqslant 2\|\lambda(x(t)-x^*(t))+\dot{x}(t)\|^2+4(F(x(t))-F^*).
\end{aligned}\end{equation}
By applying the previous inequalities, we have that:
\begin{equation}
\|\dot{x}(t)\|^2\leqslant 4\left(1+\frac{1}{\sqrt{2}-1}\right)\mathcal{E}(t).
\end{equation}
The bound on the energy given in \eqref{eq:control_nrj_hbf} lead us to the conclusion:
\begin{equation}
\|\dot{x}(t)\|=\mathcal{O}\left(e^{-\left(1-\frac{\sqrt{2}}{2}\right)t}\right).
\end{equation}\\
\subsection{Proof of Lemma \ref{lem:right-diff}}
\label{sec:proof_right_diff}

Let $\phi'$ denote the derivative of $\phi$ when it is well defined. According to \cite{young1914note}, the function $\phi$ is differentiable except at a countable set of points. This implies that there exists $(t_i)_{i\in\llbracket 1,N\rrbracket}$ and $N\in\N^*\cup\{+\infty\}$ such that for any $i\in \llbracket 0,N-1\rrbracket$ and $t\in (t_i,t_{i+1})$, $\phi'(t)$ is well defined and equal to $\phi_+(t)$. We suppose that the sequence is ordered such that $t_0<t_i< t_{i+1}$ for any $i$ and that $t_N=+\infty$ when $N\neq+\infty$. \\
Suppose that $t\in (t_0,t_1)$. 
\begin{itemize}
\item If $\phi$ is differentiable at $t_0$, then $\phi$ is differentiable on the interval $[t_0,t_1)$ and $\phi'=\phi_+$ in this interval. Consequently inequality \eqref{eq:ineq_phi+} ensures that,
\begin{equation*}
\phi(t)\leqslant \phi(t_0)+\int_{t_0}^t\psi(u)du.
\end{equation*}
\item If $\phi$ is not differentiable at $t_0$, then inequality \eqref{eq:ineq_phi+} guarantees that for $h>0$ sufficiently small,
\begin{equation*}
\phi(t_0+h)\leqslant \phi(t_0)+h\psi(t_0).
\end{equation*}
Then, the previous discussion allows us to say that $\phi$ is differentiable on $[t_0+h, t_1)$. As a consequence, we can say that there exists $H\in(0,t-t_0)$ such that for any $h\in(0,H)$:
\begin{equation*}
\phi(t)\leqslant \phi(t_0+h)+\int_{t_0+h}^t\psi(u)du\leqslant \phi(t_0)+\int_{t_0}^t\psi(u)du+\int_{t_0}^{t_0+h}\left(\psi(t_0)-\psi(u)\right)du.
\end{equation*}
As this inequality is valid for any $h\in(0,H)$, we finally get the wanted inequality \eqref{eq:maj_semidif}.
\end{itemize}
We now suppose that $t=t_1$. We just proved that \eqref{eq:maj_semidif} is true for all $t\in(t_0,t_1)$. Therefore, for all $t\in (t_0,t_1)$,
\begin{equation*}
\phi(t)\leqslant \phi(t_0)+\int_{t_0}^{t_1}\psi(u)du,
\end{equation*}
and as $\phi$ is continuous we get the same inequality at $t=t_1$.\\
By using the same arguments, we can prove that \eqref{eq:maj_semidif} is valid for any $t>t_1$. Indeed, if $t>t_1$, then it means that $t\in (t_i,t_{i+1})$ or that $t=t_i$ for some $i\in\llbracket 1,N\rrbracket$. In both cases, we get the wanted inequality by applying the above reasonings to the consecutive intervals $(t_j,t_{j+1})$ for $0\leqslant j\leqslant i$.\\

\subsection{Proof of Lemma \ref{lem:tech1}}
\label{sec:proof_lemma_tech1}
Let $n\in\N^*$. By rewriting \begin{equation*}x_n-x_n^*=\frac{1}{2}\left((x_n-x_{n-1})+(x_{n-1}-x_{n-1}^*)+(x_{n-1}^*-x_n^*)+(x_n-x_n^*)\right),\end{equation*}we get that:
\begin{equation*}
\langle x_n-x_n^*,x_n-x_{n-1}\rangle=\frac{1}{2}\delta_n+\frac{1}{2}\langle(x_{n-1}-x_{n-1}^*)+(x_{n-1}^*-x_n^*)+(x_n-x_n^*),x_n-x_{n-1}\rangle.
\end{equation*}
Noticing that $x_n-x_{n-1}=(x_n-x_n^*)+(x_n^*-x_{n-1}^*)+(x_{n-1}^*-x_{n-1})$ leads to:
\begin{equation*}
\begin{aligned}
2\langle x_n-x_n^*,x_n-x_{n-1}\rangle&=\delta_n+\langle x_{n-1}-x_{n-1}^*,x_n-x_n^*\rangle +\langle x_{n-1}-x_{n-1}^*,x_n^*-x_{n-1}^*\rangle\\
&-h_{n-1}-\langle x_n^*-x_{n-1}^*,x_n-x_n^*\rangle+\langle x_{n-1}-x_{n-1}^*,x_n^*-x_{n-1}^*\rangle\\
&-\gamma_n^*+\langle x_{n-1}-x_{n-1}^*,x_n^*-x_{n-1}^*\rangle-\langle x_{n-1}-x_{n-1}^*,x_n-x_n^*\rangle+h_n\\
&=h_n-h_{n-1}+\delta_n-\gamma_n^*+2\langle x_{n-1}-x_{n-1}^*,x_n^*-x_{n-1}^*\rangle.
\end{aligned}
\end{equation*}
The second claim is proved using the same approach. We rewrite \begin{equation*}x_{n-1}-x_{n-1}^*=\frac{1}{2}\left((x_{n-1}-x_n)+(x_n-x_n^*)+(x_n^*-x_{n-1}^*)+(x_{n-1}^*-x_{n-1})\right),\end{equation*}and consequently:
\begin{equation*}
2\langle x_{n-1}-x_{n-1}^*,x_n-x_{n-1}\rangle=-\delta_n+\langle (x_n-x_n^*)+(x_n^*-x_{n-1}^*)+(x_{n-1}^*-x_{n-1}),x_n-x_{n-1}\rangle.
\end{equation*}
By applying the same rewriting of $x_n-x_{n-1}$, simple calculations give that:
\begin{equation*}
\langle x_{n-1}-x_{n-1}^*,x_n-x_{n-1}\rangle=\frac{1}{2}(h_n-h_{n-1}-\delta_n+\gamma_n^*)+\langle x_{n}-x_{n}^*,x_n^*-x_{n-1}^*\rangle.
\end{equation*}

\subsection{Proof of Lemma \ref{lem:tech2}}
\label{sec:proof_lemma_tech2}
The first claim is straight forward as Lemma 3.1 of \cite{chambolle2015convergence} ensures that:
\begin{equation*}
F(x_{n+1})-F(x_n)\leqslant\frac{L}{2}\left(\|y_n-x_n\|^2-\|x_{n+1}-x_n\|^2\right).
\end{equation*}
By writing $y_n=x_n+\alpha_n(x_n-x_{n-1})$ and $\frac{2}{L}(F(x_{n+1})-F(x_n))=w_{n+1}-w_n$, we can conclude.

By applying Lemma 3.1 of \cite{chambolle2015convergence} to an other couple of points, we get that:
\begin{equation*}
F(x_{n+1})-F^*\leqslant \frac{L}{2}\left(\|y_n-x_n^*\|^2-\|x_{n+1}-x_n^*\|^2\right).
\end{equation*}
It follows that:
\begin{equation*}
\begin{aligned}
w_{n+1}&\leqslant\|x_n+\alpha_n(x_n-x_{n-1})-x_n^*\|^2-\|(x_{n+1}-x_{n+1}^*)+(x_{n+1}^*-x_n^*)\|^2\\
&\leqslant h_n+\alpha_n^2\delta_n-h_{n+1}-\gamma_{n+1}^*+2\alpha_n\langle x_n-x_n^*,x_n-x_{n-1}\rangle\\&-2\langle x_{n+1}-x_{n+1}^*,x_{n+1}^*-x_n^*\rangle.
\end{aligned}
\end{equation*}
Recall that the first claim of Lemma \ref{lem:tech1} ensures that:
\begin{equation*}
\langle x_n-x_n^*,x_n-x_{n-1}\rangle=\frac{1}{2}(h_n-h_{n-1}+\delta_n-\gamma_n^*)+\langle x_{n-1}-x_{n-1}^*,x_n^*-x_{n-1}^*\rangle,
\end{equation*}
we can deduce that:
\begin{equation*}
\begin{aligned}
w_{n+1}&\leqslant (1+\alpha_n)h_n+(\alpha_n^2+\alpha_n)\delta_n-\alpha_n h_{n-1}-h_{n+1}-\gamma_{n+1}^*-\alpha_n\gamma_n^*\\
&+2\alpha_n\langle x_{n-1}-x_{n-1}^*,x_n^*-x_{n-1}^*\rangle-2\langle x_{n+1}-x_{n+1}^*,x_{n+1}^*-x_n^*\rangle.
\end{aligned}
\end{equation*}

\subsection{Proof of Lemma \ref{lem:VFISTA_nrj}}
\label{sec:proof_lem_nrj}
Let $(x_n)_{n\in\N}$ be the sequence provided by \eqref{eq:V-FISTA} and $s=\frac{1}{L}$. We can write the Lyapunov energy $(\mathcal{E}_n)_{n\in\N}$ in the following way:
\begin{equation*}
\mathcal{E}_n=w_n+(1-\lambda)\delta_n+\lambda (h_n-h_{n-1})+\lambda^2h_{n-1}+\lambda\gamma_n^*+2\lambda\langle x_n-x_n^*,x_n^*-x_{n-1}^*\rangle.
\end{equation*}
Since $\langle x_n-x_n^*,x_n^*-x_{n-1}^*\rangle\geqslant0$, we can write that:
\begin{equation*}
\mathcal{E}_n\geqslant w_n+\lambda (h_n-h_{n-1}),
\end{equation*}
which leads to the final result.

\subsection{Proof of Lemma \ref{lem:VFISTA_diff}}
\label{sec:proof_lem_diff_vf2}

Let $(x_n)_{n\in\N}$ be the sequence provided by \eqref{eq:V-FISTA}. By using the expression \eqref{eq:rewriting_lyap2} of $\mathcal{E}_n$, we get that:
\begin{equation}
\begin{aligned}
\mathcal{E}_{n+1}-\mathcal{E}_n=&~w_{n+1}-w_n+\lambda\left(\alpha+\lambda\alpha+(1-\alpha)^2\right)(h_{n+1}-h_n)+\alpha(1+\lambda)\delta_{n+1}\\&-\lambda\alpha(h_n-h_{n-1})-\alpha(1+\lambda)\delta_n-\lambda\alpha\gamma_{n+1}^*+\lambda\alpha\gamma_n^*\\&+2\lambda\alpha\langle x_n-x_n^*,x_{n+1}^*-x_n^*\rangle-2\lambda\alpha\langle x_{n-1}-x_{n-1}^*,x_n^*-x_{n-1}^*\rangle.
\end{aligned}
\end{equation}
The first claim of Lemma \ref{lem:tech2} combined to the inequality $-\lambda\alpha\gamma_{n+1}^*+2\lambda\alpha\langle x_n-x_n^*,x_{n+1}^*-x_n^*\rangle\leqslant0$ lead to:
\begin{equation}
\begin{aligned}
\mathcal{E}_{n+1}-\mathcal{E}_n\leqslant&~\lambda\left(\alpha+\lambda\alpha+(1-\alpha)^2\right)(h_{n+1}-h_n)+(\alpha+\lambda\alpha-1)\delta_{n+1}\\&-\lambda\alpha(h_n-h_{n-1})-\alpha(1+\lambda-\alpha)\delta_n+\lambda\alpha\gamma_n^*\\&-2\lambda\alpha\langle x_{n-1}-x_{n-1}^*,x_n^*-x_{n-1}^*\rangle.
\end{aligned}
\end{equation}
According to the second claim of Lemma \ref{lem:tech2}, we have
\begin{equation}
\begin{aligned}
0\leqslant &~-\lambda w_{n+1}+\lambda(\alpha+\alpha^2)\delta_n+\lambda\alpha(h_n-h_{n-1})-\lambda(h_{n+1}-h_n)\\
&-\lambda\gamma_{n+1}^*-\lambda\alpha\gamma_n^*+2\lambda\alpha\langle x_{n-1}-x_{n-1}^*,x_n^*-x_{n-1}^*\rangle\\&-2\lambda\langle x_{n+1}-x_{n+1}^*,x_{n+1}^*-x_n^*\rangle,
\end{aligned}
\end{equation}
and as $\langle x_{n+1}-x_{n+1}^*,x_{n+1}^*-x_n^*\rangle\geqslant0$, it follows that
\begin{equation}
\begin{aligned}
\mathcal{E}_{n+1}-\mathcal{E}_n\leqslant&~-\lambda w_{n+1}+\lambda\left(\alpha+\lambda\alpha+(1-\alpha)^2-1\right)(h_{n+1}-h_n)\\&+(\alpha+\lambda\alpha-1)\delta_{n+1}-\alpha(1-\alpha-\lambda\alpha)\delta_n.
\end{aligned}
\end{equation}
By developing $\alpha+\lambda\alpha+(1-\alpha)^2-1$ we get to the conclusion.

\section*{Acknowledgements}

JFA acknowledges support of the  EU  Horizon  2020 research and innovation program under the Marie Sk\l odowska-Curie NoMADS grant agreement No777826, and PEPR PDE-AI. HL acknowledges the financial support of the Ministry of Education, University and Research (grant ML4IP R205T7J2KP). This work was supported by the ANR MICROBLIND (grant ANR-21-CE48-0008) and the ANR Masdol (grant ANR-PRC-CE23). 

 \bibliographystyle{abbrv}
\bibliography{ref.bib}

\begin{thebibliography}{10}

\bibitem{alamo2022restart}
T.~Alamo, P.~Krupa, and D.~Limon.
\newblock Restart of accelerated first-order methods with linear convergence under a quadratic functional growth condition.
\newblock {\em IEEE Transactions on Automatic Control}, 68(1):612--619, 2022.

\bibitem{apidopoulos2022convergence}
V.~Apidopoulos, N.~Ginatta, and S.~Villa.
\newblock Convergence rates for the {H}eavy-{B}all continuous dynamics for non-convex optimization, under {P}olyak--{{\L}}ojasiewicz condition.
\newblock {\em Journal of Global Optimization}, pages 1--27, 2022.

\bibitem{attouch2000heavy}
H.~Attouch, X.~Goudou, and P.~Redont.
\newblock The {H}eavy-{B}all with friction method, {I}. the continuous dynamical system: global exploration of the local minima of a real-valued function by asymptotic analysis of a dissipative dynamical system.
\newblock {\em Communications in Contemporary Mathematics}, 2(01):1--34, 2000.

\bibitem{attouch2016rate}
H.~Attouch and J.~Peypouquet.
\newblock The rate of convergence of nesterov's accelerated forward-backward method is actually faster than $1/k^2$.
\newblock {\em SIAM Journal on Optimization}, 26(3):1824--1834, 2016.

\bibitem{aujol2023parameter}
J.-F. Aujol, L.~Calatroni, C.~Dossal, H.~Labarri{\`e}re, and A.~Rondepierre.
\newblock Parameter-free {FISTA} by adaptive restart and backtracking.
\newblock {\em arXiv preprint arXiv:2307.14323}, 2023.

\bibitem{aujol2021restart}
J.-F. Aujol, C.~Dossal, H.~Labarri{\`e}re, and A.~Rondepierre.
\newblock {{FISTA} restart using an automatic estimation of the growth parameter}.
\newblock {\em Hal Preprint 03153525}, May 2022.

\bibitem{aujol2024strong}
J.-F. Aujol, C.~Dossal, H.~Labarri{\`e}re, and A.~Rondepierre.
\newblock Strong convergence of {FISTA} {I}terates under {H}\"olderian and {Q}uadratic {G}rowth {C}onditions.
\newblock {\em arXiv preprint arXiv:2407.17063}, 2024.

\bibitem{aujol2022convergenceQC}
J.-F. Aujol, C.~Dossal, and A.~Rondepierre.
\newblock Convergence rates of the {H}eavy {B}all method for quasi-strongly convex optimization.
\newblock {\em SIAM Journal on Optimization}, 32(3):1817--1842, 2022.

\bibitem{aujol2022convergenceL}
J.-F. Aujol, C.~Dossal, and A.~Rondepierre.
\newblock Convergence rates of the {H}eavy-{B}all method under the {{\L}}ojasiewicz property.
\newblock {\em Mathematical Programming}, pages 1--60, 2022.

\bibitem{Aujol2023}
J.-F. Aujol, C.~Dossal, and A.~Rondepierre.
\newblock {FISTA} is an automatic geometrically optimized algorithm for strongly convex functions.
\newblock {\em Mathematical Programming}, 204(1-2), 2024.

\bibitem{beck2017first}
A.~Beck.
\newblock {\em First-order methods in optimization}.
\newblock SIAM, 2017.

\bibitem{beck2009fast}
A.~Beck and M.~Teboulle.
\newblock A fast iterative shrinkage-thresholding algorithm for linear inverse problems.
\newblock {\em SIAM journal on imaging sciences}, 2(1):183--202, 2009.

\bibitem{begout2015damped}
P.~B{\'e}gout, J.~Bolte, and M.~A. Jendoubi.
\newblock On damped second-order gradient systems.
\newblock {\em Journal of Differential Equations}, 259(7):3115--3143, 2015.

\bibitem{Bolte2007loja}
J.~Bolte, A.~Daniilidis, A.~Lewis, and M.~Shiota.
\newblock Clarke subgradients of stratifiable functions.
\newblock {\em SIAM Journal on Optimization}, 18(2):556--572, 2007.

\bibitem{bolte2017error}
J.~Bolte, T.~P. Nguyen, J.~Peypouquet, and B.~W. Suter.
\newblock From error bounds to the complexity of first-order descent methods for convex functions.
\newblock {\em Mathematical Programming}, 165:471--507, 2017.

\bibitem{bonnans1998sensitivity}
J.~F. Bonnans, R.~Cominetti, and A.~Shapiro.
\newblock Sensitivity analysis of optimization problems under second order regular constraints.
\newblock {\em Mathematics of Operations Research}, 23(4):806--831, 1998.

\bibitem{bubeck2015geometric}
S.~Bubeck, Y.~T. Lee, and M.~Singh.
\newblock A geometric alternative to {N}esterov's accelerated gradient descent.
\newblock {\em arXiv preprint arXiv:1506.08187}, 2015.

\bibitem{chambolle2015convergence}
A.~Chambolle and C.~Dossal.
\newblock On the convergence of the iterates of the “fast iterative shrinkage/thresholding algorithm”.
\newblock {\em Journal of Optimization theory and Applications}, 166(3):968--982, 2015.

\bibitem{chen2017geometric}
S.~Chen, S.~Ma, and W.~Liu.
\newblock Geometric descent method for convex composite minimization.
\newblock {\em Advances in Neural Information Processing Systems}, 30, 2017.

\bibitem{Drori2014}
Y.~Drori and M.~Teboulle.
\newblock Performance of first-order methods for smooth convex minimization: a novel approach.
\newblock {\em Mathematical Programming}, 145(1):451--482, June 2014.

\bibitem{fercoq2019adaptive}
O.~Fercoq and Z.~Qu.
\newblock Adaptive restart of accelerated gradient methods under local quadratic growth condition.
\newblock {\em IMA Journal of Numerical Analysis}, 39(4):2069--2095, 2019.

\bibitem{garrigos2022convergence}
G.~Garrigos, L.~Rosasco, and S.~Villa.
\newblock Convergence of the forward-backward algorithm: beyond the worst-case with the help of geometry.
\newblock {\em Mathematical Programming}, pages 1--60, 2022.

\bibitem{ghadimi2015global}
E.~Ghadimi, H.~R. Feyzmahdavian, and M.~Johansson.
\newblock Global convergence of the {H}eavy-{B}all method for convex optimization.
\newblock In {\em 2015 European control conference (ECC)}, pages 310--315. IEEE, 2015.

\bibitem{giselsson2014monotonicity}
P.~Giselsson and S.~Boyd.
\newblock Monotonicity and restart in fast gradient methods.
\newblock In {\em 53rd IEEE Conference on Decision and Control}, pages 5058--5063. IEEE, 2014.

\bibitem{hiriart1982points}
J.-B. Hiriart-Urruty.
\newblock At what points is the projection mapping differentiable?
\newblock {\em The American Mathematical Monthly}, 89(7):456--458, 1982.

\bibitem{kassing2024polyak}
S.~Kassing and S.~Weissmann.
\newblock Polyak's heavy ball method achieves accelerated local rate of convergence under polyak-lojasiewicz inequality.
\newblock {\em arXiv preprint arXiv:2410.16849}, 2024.

\bibitem{Loja63}
S.~{\L}ojasiewicz.
\newblock Une propri{\'e}t{\'e} topologique des sous-ensembles analytiques r{\'e}els.
\newblock In {\em Les \'{E}quations aux {D}{\'e}riv{\'e}es {P}artielles ({P}aris, 1962)}, pages 87--89. {\'E}ditions du Centre National de la Recherche Scientifique, Paris, 1963.

\bibitem{Loja93}
S.~{\L}ojasiewicz.
\newblock Sur la g{\'e}om{\'e}trie semi- et sous-analytique.
\newblock {\em Annales de l'Institut Fourier. Universit{\'e} de Grenoble}, 43(5):1575--1595, 1993.

\bibitem{necoara2019linear}
I.~Necoara, Y.~Nesterov, and F.~Glineur.
\newblock Linear convergence of first order methods for non-strongly convex optimization.
\newblock {\em Mathematical Programming}, 175(1):69--107, 2019.

\bibitem{nesterov1983method}
Y.~Nesterov.
\newblock A method of solving a convex programming problem with convergence rate o (1/k$^2$).
\newblock In {\em Sov. Math. Dokl}, volume~27, 1983.

\bibitem{nesterov2003introductory}
Y.~Nesterov.
\newblock {\em Introductory lectures on convex optimization: A basic course}, volume~87.
\newblock Springer Science \& Business Media, 2003.

\bibitem{nesterov2013gradient}
Y.~Nesterov.
\newblock Gradient methods for minimizing composite functions.
\newblock {\em Mathematical programming}, 140(1):125--161, 2013.

\bibitem{o2015adaptive}
B.~O’donoghue and E.~Candes.
\newblock Adaptive restart for accelerated gradient schemes.
\newblock {\em Foundations of computational mathematics}, 15(3):715--732, 2015.

\bibitem{polyakshch}
B.~Polyak and P.~Shcherbakov.
\newblock Lyapunov functions: An optimization theory perspective.
\newblock {\em IFAC-PapersOnLine}, 50(1):7456--7461, 2017.
\newblock 20th IFAC World Congress.

\bibitem{polyak1964some}
B.~T. Polyak.
\newblock Some methods of speeding up the convergence of iteration methods.
\newblock {\em Ussr computational mathematics and mathematical physics}, 4(5):1--17, 1964.

\bibitem{renegar2022simple}
J.~Renegar and B.~Grimmer.
\newblock A simple nearly optimal restart scheme for speeding up first-order methods.
\newblock {\em Foundations of computational mathematics}, 22(1):211--256, 2022.

\bibitem{shapiro2016differentiability}
A.~Shapiro.
\newblock Differentiability properties of metric projections onto convex sets.
\newblock {\em Journal of Optimization Theory and Applications}, 169(3):953--964, 2016.

\bibitem{siegel2019accelerated}
J.~W. Siegel.
\newblock Accelerated first-order methods: Differential equations and {L}yapunov functions.
\newblock {\em arXiv preprint arXiv:1903.05671}, 2019.

\bibitem{taylor2022optimal}
A.~Taylor and Y.~Drori.
\newblock An optimal gradient method for smooth strongly convex minimization.
\newblock {\em Mathematical Programming}, pages 1--38, 2022.

\bibitem{van2017fastest}
B.~Van~Scoy, R.~A. Freeman, and K.~M. Lynch.
\newblock The fastest known globally convergent first-order method for minimizing strongly convex functions.
\newblock {\em IEEE Control Systems Letters}, 2(1):49--54, 2017.

\bibitem{young1914note}
G.~C. Young.
\newblock A note on derivates and differential coefficients.
\newblock {\em Acta mathematica}, 37(1):141--154, 1914.

\end{thebibliography}

\end{document}